\documentclass{article}
\usepackage[utf8]{inputenc}
\usepackage[margin=1in]{geometry}
\usepackage{amsmath} % assumes amsmath package installed
\usepackage{todonotes}
\usepackage{xcolor}
\usepackage{amsmath}
\usepackage{amssymb}
\usepackage{mathtools}
\usepackage{subcaption}
\usepackage{hyperref}
\hypersetup{colorlinks=true, linkcolor=black}
\usepackage{algorithm}
\usepackage{algpseudocode}
\usepackage{mathrsfs}

\usepackage{booktabs}   % For \toprule, \midrule, \bottomrule, \cmidrule
\usepackage{multirow} 

\usepackage{tikz,pgfplots}
\usepackage{color,graphicx,fancybox,fancyvrb}

\definecolor{ccolor}{RGB}{203,96,21}
\newcommand{\rw}[1]{ #1}

\definecolor{matlab1}{RGB}{0  113.9850  188.9550}
\definecolor{matlab2}{RGB}{ 216.7500   82.8750   24.9900}
\definecolor{matlab3}{RGB}{  236.8950  176.9700   31.8750}
\definecolor{matlab4}{RGB}{  125.9700   46.9200  141.7800}

\definecolor{matlab5}{RGB}{    118.8300  171.8700   47.9400}

\newcommand{\R}{\mathbb{R}}

\newcommand{\N}{\mathbb{N}}

\newcommand{\s}{\mathcal{S}}

\newcommand{\0}{\mathbf{0}}

\newcommand{\psd}{\mathbb{S}}

\newcommand{\vs}{\mathcal{V}}
\newcommand{\es}{\mathcal{E}}
\newcommand{\gs}{\mathcal{G}}

\DeclarePairedDelimiter{\abs}{\lvert}{\rvert}
\DeclarePairedDelimiter{\norm}{\lVert}{\rVert}

\DeclarePairedDelimiterX{\inp}[2]{\langle}{\rangle}{#1, #2}

\DeclareMathOperator*{\argmin}{arg\!\,min}

\definecolor{cn}{RGB}{93,147,191}           %0
\definecolor{cz}{RGB}{233,  72, 73}      %1
\definecolor{cp}{RGB}{113, 191, 110}    %2
\definecolor{cf}{RGB}{157,0,255}    %2

\newtheorem{thm}{Theorem}[section]
\newtheorem{lem}[thm]{Lemma}

\newtheorem{prob}{Problem}
\newtheorem{assum}{Assumption}
\newtheorem{cor}{Corollary}

\newtheorem{remark}{Remark} 

\newcommand{\sml}{\mathcal{S}_{m, L}}

\usepackage{tikz,pgfplots,datatool}
\usepackage[tikz]{bclogo}
\newcommand{\skipthis}[1]{}

\usepackage{arydshln}

\newcommand{\hdl}{\\ \hdashline}
\newcommand{\mat}[2]{\left(\begin{array}{#1}#2\end{array}\right)}
\definecolor{yel}{rgb}{1,.9,.6}
\definecolor{lmag}{rgb}{0,0.82,1}

\definecolor{mag}{rgb}{0,0.7,.9}
\definecolor{magenta}{rgb}{1,0,1}
\definecolor{lgrey}{rgb}{0.85,0.85,0.85}
\definecolor{lgr}{rgb}{.9,0.9,0.9}
\definecolor{lgreen}{rgb}{0,0.92,0.7}
\colorlet{lred}{red!20!white}
\definecolor{lblue}{rgb}{1,0.95,0.7}

\definecolor{sblue}{rgb}{.1,.2,.5}
\colorlet{lsblue}{sblue!20!white}
\definecolor{sred}{rgb}{.85,.25,0}
\colorlet{lsred}{sred!20!white}
\definecolor{lye}{rgb}{1,0.95,0.7}

\usepgflibrary{fpu}
\usetikzlibrary{positioning,calc,arrows,shapes,shadows,matrix,backgrounds,external}
\tikzsetexternalprefix{figures/}
\tikzset{
auto,
poi/.style={
minimum size=0,
inner sep=0
},
sys/.style 2 args={
rectangle,
draw,
rounded corners,
drop shadow,
draw=black,
top color=black!20,bottom color=black!0,
minimum height=#2,
minimum width=#1,
inner sep=\dn},
syse/.style 2 args={
rectangle,
draw=none,
rounded corners,
minimum height=#2,
minimum width=#1,
inner sep=\dn},
nod/.style={
circle,
draw,
fill=white,
minimum size=5ex
},
sum/.style={circle,draw,draw=black,inner sep=0mm,minimum size=2mm,drop shadow,fill=white,
draw=black!100,top color=black!20,bottom color=black!0},
sume/.style={circle,draw=none,inner sep=0mm,minimum size=2mm,drop shadow,fill=white,
draw=black!100,top color=black!20,bottom color=black!0},
jun/.style={circle,draw,draw=black,inner sep=0mm,minimum size=0mm},
>={latex},
every path/.style={rounded corners},
lin/.style={color=black,draw,->},
}

\def\dn{1ex}

\tikzstyle{sy0}=[sys={0*\dn}{0*\dn}]
\tikzstyle{sy1}=[sys={12*\dn}{8*\dn}]
\tikzstyle{sy2}=[sys={8*\dn}{6*\dn}]
\tikzstyle{sy3}=[sys={5*\dn}{5*\dn}]
\tikzstyle{sy0}=[sys={0*\dn}{0*\dn}]
\tikzstyle{sye1}=[syse={12*\dn}{8*\dn}]
\tikzstyle{sye2}=[syse={8*\dn}{6*\dn}]
\tikzstyle{sye3}=[syse={5*\dn}{5*\dn}]

\tikzstyle{sy4}=[sys={12*\dn}{6*\dn}]
\tikzstyle{sye4}=[syse={12*\dn}{6*\dn}]

\newcommand{\te}[1]{\text{\ \ #1\ \ }}

\newcommand{\Acl}{{\cal A}}
\newcommand{\Bcl}{{\cal B}}
\newcommand{\Ccl}{{\cal C}}
\newcommand{\Dcl}{{\cal D}}

\newcommand{\hl}{\\\hline}

\newcommand{\qed}{\hfill$\square$}
\newcommand{\citep}[1]{\cite{#1}}
\renewcommand\footnotemark{}

\title{Analysis and Synthesis of Switched Optimization Algorithms}

\author{Jared Miller$^1$, Carsten Scherer$^1$, Fabian Jakob$^2$, Andrea Iannelli$^2$
\thanks{ J. Miller and C. Scherer are with the Chair of Mathematical Systems Theory, Department of Mathematics,  University of Stuttgart, Stuttgart, Germany 
  (\{jared.miller, carsten.scherer\}@imng.uni-stuttgart.de)
  }  
  \thanks{{ F. Jakob and A. Iannelli are with the Institute for Systems Theory and Automatic Control, University of Stuttgart, Germany
(e-mail: \{fabian.jakob, andrea.iannelli\}@ist.uni-stuttgart.de).}}
  \thanks{J. Miller and C. Scherer are funded by Deutsche Forschungsgemeinschaft (DFG, German Research Foundation) under Germany's Excellence Strategy - EXC 2075 – 390740016. J. Miller and C. Scherer acknowledge the support by the Stuttgart Center for Simulation Science (SimTech). F. Jakob acknowledges the support of the International Max Planck Research School for Intelligent Systems (IMPRS-IS).}}

\begin{document}

\maketitle
% \thispagestyle{empty}
% \pagestyle{empty}

%%%%%%%%%%%%%%%%%%%%%%%%%%%%%%%%%%%%%%%%%%%%%%%%%%%%%%%%%%%%%%

\begin{abstract}
\label{sec:abstract}
Deployment of optimization algorithms over \rw{communication networks} face challenges associated with time delays and corruptions. 
% One particular instance is the presence of time-varying delays  arising  from factors such as packet drops. 
Fixed time delays can destabilize popular \rw{gradient-based }algorithms, and this degradation is exacerbated by time-varying delays \rw{that may arise from packet drops}.  This work concentrates on the analysis and \rw{synthesis} of discrete-time optimization algorithms with certified exponential convergence rates that are robust against switched network dynamics between the optimizer and the gradient oracle. 
% These optimization algorithms are implemented by switch-scheduled output feedback controllers.
% and their convergence that is robust 
% Optimization algorithm synthesis is accomplished by computing a set of switching-dependent output-feedback controllers, that,  together ensures the closed-loop system satisfies a passivity requirement. 
Analysis is accomplished by \rw{solving linear matrix inequalities} under bisection in the exponential convergence rate, searching over  Zames-Falb filter coefficients that can certify convergence.  Synthesis is performed by alternating between a search over filter coefficient for a fixed controller, and a search over controllers for a fixed \rw{filter. Effectiveness is demonstrated by the  synthesis of convergent optimization algorithms over networks with time-varying delays,  and networks with unstable channel dynamics.}

 % solving for Zames-Falb filter coefficients to certify a minimal convergence rate, and Synthesis is performed by alternation 

% The analysis task is accomplished by solving

% synthesis task is accomplished by solving linear matrix inequalities in an alternating manner.

% by solving semidefinite programs associated with ensuring that a gain-scheduled linear-parameter-varying system 

% The presence of time-delays may destabilize optimization algorithms that are cerr
% \urg{Abstract goes here}

\end{abstract}


%% Use \section commands to start a section
\section{Introduction}
\label{sec:introduction}

Optimization algorithms increasingly operate within networked environments, where information exchange is rarely perfect. Communication links introduce delays, packet drops, and other distortions that alter how optimization algorithms interact with \rw{gradients}. These effects can have significant consequences -- slowing convergence, degrading performance, or even destabilizing the algorithm dynamics.

% To address these issues, recent works have investigated the design of optimization methods that are robust to uncertainty in the oracle interaction. 
In the strongly convex setting, robust control techniques have been used to analyze and synthesize optimization algorithms.
Introduced by \cite{lessard2016analysis}, the essential idea is to treat a first-order algorithm as \rw{the interconnection of a linear system with a static gradient nonlinearity}, and to frame the convergence as an absolute stability problem \rw{of the resultant Lur'e system \cite{lur1944theory}}. 
\rw{Algorithms should be designed such that they are applicable to a class of objective functions. Input-output relations satisfied by sequences arising from gradients inside the function class} 
 are described using Integral Quadratic Constraints (IQCs)
\citep{Megretski1997,BoczarLessardRecht}, and the \rw{analysis and} design of robust algorithms can subsequently be cast as a \rw{robust control} problem \citep{scherer2021convex}. This methodology \rw{also} enables \rw{a quantification of}  robustness to \rw{noise in the gradients as measured by}   $H_2$ performance \citep{michalowski2021,jovanovic2021robustness}, or to account for adversarial communication dynamics by considering them within an extremum control setup \citep{holicki2021algorithm}. 

Most of this system-theoretic approach has focused on algorithm and communication dynamics \rw{that are each time-invariant} \citep{scherer2021convex, scherer2023optimization}.
%, leaving network-induced phenomena such as \emph{time-varying delays} largely unaddressed. 
Even though recently this approach was extended to time-varying objectives and algorithms in \citep{jakob2025linear}, important network-induced phenomena such as \emph{time-varying delays} have not been specifically addressed yet.
%in automated algorithm synthesis. 
When such effects are taken
into account, the resulting dynamics are more accurately
captured by switched  \rw{systems} \citep{wen2008switched, conte2020modeling}.
% As such, methodologies from control switched 
% for which rich synthesis methodologies already exist \citep{wen2008switched, liu2011delay}. 
% In particular, gain-scheduled or LPV control frameworks allow for controllers that dynamically adapt to network delay \citep{de2012gain}. Incorporating similar ideas in the optimization framework could enable algorithms to adapt to variations in delay, rather than merely tolerating it.
Motivated by this perspective, this work extends the system-theoretic framework to include \emph{switched optimization algorithms}.
% Switched optimization algorithms arise in the setting of time-varying delays , and \citep{wen2008switched, conte2020modeling,conte2021disturbance}
% are also found when applying a time-independent algorithm (e.g. standard gradient descent) over a time-varying network. 
The design of optimization algorithms in a switched setting
has mostly focused on analysis and synthesis of distributed algorithms that are robust against time-varying information delays in the form of stale (time-delayed) gradients \citep{ramaswamy2021optimization,  doostmohammadian2021consensus}. 
Delay-scheduled algorithms, in analogy to delay-scheduled controllers \citep{briat2009delay}, have not yet been addressed.
%A switched optimization algorithm could also be delay-scheduled \citep{briat2009delay}, in which the algorithm parameters would change based on the measured round-trip delay from the remote gradient.

% has mostly focused on stepsize selection of first and second order gradient descent methods in a asynchronous distributed setting that are robust 
% Because time-varying network delays can be described by switched systems , we synthesize delay-scheduled \citep{briat2009delay} optimization algorithms by using tools from robust control. 
% Given that 
% /(e.g. measured delay). 
% We discuss structural constraints that an algorithm must fulfill to achieve convergence in a switched setting, and decompose the algorithm into an internal model for reference tracking and a free-to-synthesize subsystem. 
% We use parameter-affine \citep{apkarian1995self}/path-dependent \citep{Ahmadi_2014} storage functions 

This work performs analysis of switched optimization algorithms by solving mode-dependent \citep{apkarian1995self} linear matrix inequalities (LMIs) \citep{SCHERER2023robustexp}. 
The switched optimization synthesis problem is posed as a robustness requirement under a regulation constraint \citep{francis1976internal, stoorvogel2000performance}, \rw{which is } compatible under mode changes \citep{korouglu2008lpv,  conte2021disturbance}. \rw{Obedience of the regulation constraint is ensured through introduction of a network-dependent internal model.}
% . A per-mode internal model of an integrator is sufficient to meet this 
% We use regulation arguments \cite{francis1976internal, swit} 
Given fixed \rw{Zames-Falb filters from the analysis procedure}, we perform switched optimization algorithm synthesis by applying \rw{controller design} techniques for LPV systems \citep{de2012gain}. 
The convergence rates and empirical performance of the resulting algorithms are demonstrated on networks with unstable channel dynamics, and networks with time-varying  delays.

% We note that for delay-scheduling in continuous time include \citep{briat2009delay, briat2009cal} for delay-scheduled controller synthesis with prescribed $H_\infty$ norms.

% A controller synthesis procedure based on the polytopic LPV framework and IQC-synthesis \citep{scherer2023optimization} is proposed, performing an alternating search over controller parameters and . 

% \input{sections/acronym}

The contributions of this work are:
\begin{itemize}
    \item A switched-system framework for  first-order algorithms with bilateral corrupted gradient transmission, encompassing dynamical network models such as packet dropout and bounded-rate delay.
    % \item A derivation of necessary structural regulation properties for any valid switched optimization algorithm
    \item An analysis procedure for upper-bounding worst-case convergence rates of given switched optimization algorithms.
    \item An internal-model-control-based synthesis procedure in the polytopic LPV setting for minimal worst-case linear convergence under switched network dynamics.
\end{itemize}

\section{Preliminaries}
\label{sec:preliminaries}

We review preliminaries of notation and  switched systems.

% \textbf{Notation. }

\subsection{Notation}
The set of natural numbers including 0 is $\N$. The set of natural numbers between $a$ and $b$ inclusive is $\{a, \ldots, b\}$.
%The set of natural numbers including 0 is $\N$.
%The transpose of a matrix $M \in \R^{n \times m}$ is $M^\top \in \R^{m \times n}.$  
An $\R^n$-valued signal is a map $x: \N \rightarrow \R^n$ \rw{indexed as $\{x_k\}_{k \in \N}$}. \rw{The 2-norm of a vector $g \in \R^n$ is $\norm{g}_2 = \sqrt{g^\top g}$.} 
% The one-step shift operator is $\bz$. 
\rw{A linear system $G$ is a mapping that takes an initial condition  $x_0 \in \R^n$ and an input signal $u: \N \rightarrow \R^{n_u},$ and returns an output signal $y: \N \rightarrow \R^{n_y}$. The linear system $G: (x_0, u) \rightarrow y$ may be represented by state space matrices  $(\Acl, \Bcl, \Ccl, \Dcl)$ as }
% Given an initial condition, a linear system $G$ is a mapping $(x_0, u) \rightarrow y$ with representations  of
%applying the signal $u$ to the linear system
% \begin{minipage}{0.5\linewidth}
\begin{align}
\begin{aligned}
G: \qquad \rw{\mat{c}{x_{k+1} \hl y_k}}  &= \rw{\mat{c|c}{
        \Acl & \Bcl \hl
        \Ccl & \Dcl
    } \mat{c}{x_k \hl u_k}}.
    \end{aligned} \label{eq:linear_system}
\end{align}
% \end{minipage}
% \begin{minipage}{0.5\linewidth}
%     \begin{align}
%         \rw{y = \mas{c|c}{
%         \Acl & \Bcl \hl
%         \Ccl & \Dcl} u.}
%     \end{align} \label{eq:linear_system}
% \end{minipage}
% which will be compactly expressed as 
%uniquely generates a response $(x, y)$. The notation
% will compactly express the linear system in \eqref{eq:linear_system} as a multi-valued $x_0$-dependent map $u \rightarrow y$.
 
% The series connection of linear systems $G_1: z \rightarrow w$ and $G_2: w\rightarrow p$ will be denoted as $G_1 \cdot G_2: z \rightarrow p$.
%\rw{
% The state-space matrices of a representation $G$ will be denoted as $(G)_{ss}$. 
%The transfer function of the system $G$  \rw{ in \eqref{eq:linear_system}} is $G(\bz) = \Ccl(\bz I - \Acl)^{-1} \Bcl + \Dcl.$}

% \urg{TODO: get a good citation for this.}

The symbol $e_i \in \R^n$ will denote the standard unit vector. \rw{The set of symmetric matrices of dimension $n$ is $\psd^n$.}
Positive (semi)-definitness of a symmetric matrix $M$ will be denoted as $M \succ 0$ $(M \succeq 0)$. The identity matrix of dimension $n$ is $I_n$. The zeros and all-ones matrices of dimension $m \times n$ are $0_{m \times n}$ and $1_{m \times n}$, respectively.
% Canonical transpose elements in a symmetric matrix or a quadratic form may be substituted by an asterisk $\ast$.

\subsection{Switched Systems}
% \FJ{Needs to be introduced before mentioning switched network model (?)
% }

\rw{

A switched system with $N_s$ modes, 
% with exogenous switching signal $s$ with 
$n$ states, $n_u$ inputs, and $n_y$ outputs is described by a tuple $(\Sigma, \gs)$, where $\Sigma$ collects together the $N_s$ subsystems, and $\gs$ is a directed, unweighted graph with $N_s$-vertices.} Each subsystem $\Sigma_r$ for $r \in \{1, \ldots, N_s\}$ is assumed to be described by
% Each mode $r$ of the switched system for $r \in 1, \ldots, s$ can be represented as a linear system
% The vertices $\vs$ of the graph $\gs$ are the modes $r \in 1, \ldots, s$. The edges $\es:  (r, r') \in \vs \times \vs$ of $\gs$ encode plausible switching transitions.
\begin{equation}
\label{eq:switched_system}
\begin{aligned}
    \Sigma_r: \quad & \mat{c}{x_{k+1} \hl y_k} &= \mat{c|c}{\Acl_r & \Bcl_r \hl \Ccl_r & \Dcl_r} \mat{c}{x_k \hl u_k}.
    \end{aligned}
\end{equation}
\rw{The notation $\Sigma_r$ will refer to the specific representation in \eqref{eq:switched_system} at each $r \in 1, \ldots, s$.}
% We stay in a setting where \eqref{eq:switched_system} that the tuple $(n, n_u, n_y)$ is the same between all subsystems, 
% We stay in a setting where $(N_s, n, n_u, n_y)$ are all finite. 
The vertices $\vs$ of the switching graph $\gs$ are the modes $r \in \{1, \ldots, N_s\}$. The edges $\es:  (r, r') \in \vs \times \vs$ of $\gs$ encode \rw{possible} switching transitions.  
% The switched system $(\Sigma, \Gs)$ is a mapping between an initial condition $x_0^\R^n$, a switching signal $s \in \text{Path}(\gs)$, and an input signal $u: \N \rightarrow \R^{n_u}$
% The adjacency matrix of $\es$ is $M \in\{0, 1\}^{N_s \times N_s}$, in which $M_{r, r'} = 1$ if $(r, r') \in \es$ and  $M_{r, r'} = 0$ if $(r, r')\not\in \es$.
% A switching signal $s \in \text{Path}(\gs)$
% Given a directed graph $\gs = \{\mathcal{V}, \es \}$  with vertices $\vs$ and edges $\es$, \rw{we denote the }
The set of all infinite-length paths $\{s_k\}_{k\in \N}$ in $\gs$ \rw{is denoted  as} $\text{Path}(\gs)$, in which 
\rw{$\text{Path}(\gs)$ is a subset of all maps $\N \rightarrow \{1, \ldots, N_s$\} with $(s_k, s_{k+1}) \in \es \ \ \forall k \in \N$}.
The switched system $(\Sigma, \gs)$ can therefore be interpreted as a map $(\Sigma, \gs): (x_0, s, u) \rightarrow y$. 

A trajectory of \eqref{eq:switched_system} is a tuple $(x_k, u_k, y_k, s_k)_{k \in \N}$.  \rw{This tuple is uniquely specified by an  initial condition $x_0 \in \R^n$, an input sequence $u:  \N  \rightarrow \R^{n_u}$, and a switching sequence $s \in \text{Path}(\gs)$.

}

The switched system $(\Sigma, \gs)$ is bounded (Lyapunov Stable) if there exists a $ \gamma>0$ such that  when $u=0$, for all initial points $x_0 \in \R^n$ and switching sequences $s \in \text{Path}(\gs)$, the relation $\norm{x_k}_2 \leq \gamma \norm{x_0}_2$ is obeyed.
\rw{The switched system  is asymptotically stable if when $u=0$ it holds that $\lim_{k \rightarrow \infty} x_k = 0$   for all initial points $x_0 \in \R^n$ and switching sequences $s \in \text{Path}(\gs)$.
The switched system rejects constant disturbances if it is asymptotically stable and $\lim_{k\rightarrow \infty} y_k = 0$ holds for any constant input $u$ \cite{francis1976internal}.

}

Given two linear systems $G_1: (x_1, (q, z)) \rightarrow (p, w)$ and $G_2: (x_2, (w, u)) \rightarrow (z, y)$, we denote the star product $G_1 \star G_2: ((x_1, x_2), q, u) \rightarrow (p, y)$ as the well-posed feedback interconnection of $G_1$ and $G_2$ along the common channels $w$ and $z$ \rw{as computed from given state-space representations of $G_1$ and $G_2$ \cite[Page 267-269]{zhou1998essentials}}.
% as the well-posed feedback interconnection of $G_1$ and $G_2$ along the common channels $w$ and $z$ \rw{as computed from given state-space representations of $G_1$ and $G_2$ \cite[Page 267-269]{zhou1998essentials}}.
Given two \rw{\textit{switched}} systems $(\Sigma_1, \gs): (x_1, s, (q, z)) \rightarrow (p, w)$ and $(\Sigma_2, \gs): (x_2, s, (w, u)) \rightarrow (z, y)$ \rw{sharing the same transition graph $\gs$, we denote the star product $(\Sigma_1 \star \Sigma_2, \gs): ((x_1, x_2), s, (q, u)) \rightarrow (p, y)$ as their feedback interconnection comprised of taking the star product between the per-mode linear systems as $(\Sigma_1 \star \Sigma_2)_r = \Sigma_{1, r} \star \Sigma_{2, r}$  $\forall r \in 1, \ldots, s$. The common graph $\gs$ may be omitted in the presentation of switched systems to simplify notation.}

\section{Switched Optimization Algorithms}

\label{sec:problem_formulation}

% The problem of optimizing a given function $f: \R^c \rightarrow \R$ is 

\rw{
Given scalars $0 < m < L< \infty$, we define the class $\mathcal{S}_{m,L}$ as the  set of functions $f: \R^c \rightarrow \R$ such that $f - \frac{m}{2} \norm{\cdot}_{2}^2$ and $\frac{L}{2} \norm{\cdot}_{2}^2 - f$ are both convex \cite[Definition 1.1]{rotaru2024exact}.
% We consider a class of functions 
% The set of functions we consider is the class 
% $\mathcal{S}_{m,L}$, parameterized by scalars $0 < m < L< \infty$.  The class  
In particular,  $\mathcal{S}_{m, L}$ includes the set of twice-continuously-differentiable functions $f$ that have a bounded Hessian as 
$mI \preceq \nabla^2 f(z) \preceq L I$ for all $z \in \R^c$.}
% The function $f$ is assumed to lie in a function class $\mathcal{S}_{m, L}$
% The function $f$ is $m$-strongly convex and \rw{$L$-smooth} with $0<m<L<\infty$, in the following denoted as $f \in \mathcal{S}_{m,L}$. The unique optimal solution of \eqref{eq:problem} is denoted by $z^*$, i.e., $\nabla f(z^*) = 0$. 
% The globally defined function $f: \R^c \rightarrow \R$ is assumed to be twice-continuously-differentiable with a Hessian bounded as 
% $mI \preceq \nabla^2 f(z) \leq L I$ for all $z \in \R^c$ for known constants $0 < m < L< \infty$. 
% The set of Hessian-bounded functions will be denoted as  
The optimization problem \rw{under consideration is}
\begin{align}
    z^* \in \argmin_{z \in \R^c} f(z), \label{eq:problem}
\end{align}
\rw{in which $f \in \mathcal{S}_{m, L}.$}
Because  $m > 0$, the optimum $z^*$ solving \eqref{eq:problem} is unique.
The set $\s_{m, L}^0 \subset \s_{m, L}$ is the class of functions $f \in \s_{m, L}$ such that $f(0)=0$ and $\nabla f(0) = 0$; \rw{if $f \in \s_{m, L}^0$, then $z^* = 0$ satisfies \eqref{eq:problem}.}

% The linear systems $\Sigma_r$ for each $r \in 1, \ldots, s$ each have $c$ inputs, $c$ outputs, and $n$ states, with a representation of $(\Acl_r, \Bcl_r, \Ccl_r, 0)$. Strictly properness omits proximal evaluation of the oracle $\nabla f$.

\rw{An optimization algorithm is a causal procedure that returns a sequence of iterates $\{z_k\}_{k \in \N}$ given an initial iterate $z_0 \in \R^c$ and a function $f \in \mathcal{S}_{m, L}$. The optimization algorithm is globally \textit{convergent} if for every initial iterate $z_0 \in \R^c$, the returned sequence of iterates  $\{z_k\}_{k \in \N}$ satisfies $\lim_{k \rightarrow \infty} z_k = z^*$ where $z^*$ satisfies \eqref{eq:problem}.

}
% The algorithm is convergent if $\rho < 1$.

\subsection{Algorithm Model}
% \rw{A  switched plant with $N_s < \infty$ subsystems indexed by modes may be represented by a set of linear systems $\Sigma = \{\Sigma_r\}_{r=1}^{N_s}$ and a $N_s$-vertex directed and unweighted graph $\gs$. The linear systems $\Sigma_r$ for each $r \in 1, \ldots, s$ each have $c$ inputs, $c$ outputs, and $n$ states, with a representation of $(\Acl_r, \Bcl_r, \Ccl_r, 0)$. Strictly properness omits proximal evaluation of the oracle $\nabla f$.
% }
% A trajectory of the switched system is defined with respect to a switching signal $s: \mathbb{N} \rightarrow 1, \ldots, N_s$ as $\{(s_k, x_k)\}_{k}$, in which each pair of states
% $(x_k, x_{k+1})$ obeys \eqref{eq:alg_plants} under the active subsystem $s_k$, and \rw{the signal $s$ satisfies $s \in \mathrm{Path}(\gs)$ ($ (s_k, s_{k+1}) \in \es \  \forall k \in \N)$}. 

\rw{The algorithmic interconnection of the gradient of a function  $f \in \mathcal{S}_{m, L}$ and switched system  $(\Sigma, \gs)$ with $\Dcl_{r}=0$ at all modes is
\begin{align}
   \mat{c}{x_{k+1} \hl z_k} &= \mat{c|c}{\Acl_{s_k} & \Bcl_{s_k} \hl
    \Ccl_{s_k} & 0} \mat{c}{x_{k} \hl w_k}, & w_k = \nabla f(z_k). \label{eq:alg_deployment}
\end{align}

A trajectory of the algorithmic interconnection \eqref{eq:alg_deployment} is a tuple $(x_k, w_k, z_k, s_k)_{k \in \N}$ satisfying \eqref{eq:alg_deployment} for all $k \in \N$.

% In this work, we will concentrate on cases where  $\Dcl_{s_k}=0$ (no proximal evaluation). 
% We concentrate on cases \eqref{eq:alg_deployment} in which $\Dcl_{r}=0$ for all modes $r \in 1, \ldots, s$ to avoid proximal evaluation of the oracle $f$.
% Given  $(x_0, s)$, the subsequent trajectory $(x_k, z_k, w_k, s_k)_{k \in \N}$ therefore exists and is uniquely specified (well-posed).
% A fixed point of the interconnection in \eqref{eq:alg_deployment} is a tuple $(x^*, w^*, z^*)$ such that
% \begin{align}
%    \forall (r, r') \in \es: \mat{c}{x^* \hl z^*} &= \mat{c|c}{\Acl_r & \Bcl_r\hl
%     \Ccl_r & 0 } \mat{c}{x^* \hl w}, & w^* = \nabla f(z^*). \label{eq:fixed_point_orbit}
% \end{align}
A \textit{fixed orbit} of \eqref{eq:alg_deployment} is a tuple $(\{x^*_r\}_{r=1}^{N_s}, w^*, z^*)$ such that 
\begin{align}
   \forall (r, r') \in \es: \mat{c}{x^*_{r'} \hl z^*} &= \mat{c|c}{\Acl_r & \Bcl_r\hl
    \Ccl_r & 0 } \mat{c}{x_r^* \hl w^*}, & w^* = \nabla f(z^*). \label{eq:fixed_point}
\end{align}
 % A fixed point $(x^*, w^*, z^*)$ is an instance of a fixed orbit $(\{x^*\}_{r=1}^{N_s}, w^*, z^*)$.

The interconnection in \eqref{eq:alg_deployment} is a \textit{convergent switched optimization algorithm} if there exists a fixed orbit $(\{x^*_r\}_{r=1}^{N_s}, w^*, z^*)$ such that the following properties are satisfied:
\begin{enumerate}
    \item \textbf{Optimality:} $w^* = 0$.
    \item \textbf{Zero-State:} If $z^* = 0$, then $x^*_r = 0$ for all $r \in \{1, \ldots, N_s\}$.
    \item \textbf{Attractivity:} Every trajectory $(x_k, w_k, z_k, s_k)$ satisfies
\begin{align}
    \lim_{k \rightarrow \infty} \norm{x_k - x_{s_k}^*}_2 & = 0, &  \lim_{k \rightarrow \infty} \norm{z_k - z^*}_2 &= 0, & \lim_{k \rightarrow \infty} w_k = 0. \label{eq:switched_opt_limit}
\end{align}
\end{enumerate}

% when considering all possible initial points $x_0 \in \R^n$, switching sequences $s \in \text{Path}(\mathcal{G})$, and functions $f \in \mathcal{S}_{m, L}.$

\rw{An example of a switched optimization algorithm is a gradient descent scheme with stepsize $\gamma$ under time varying delay $h_k$ as described by  $z_{k+1} = z_k-\gamma \nabla f(z_{k-h_k})$ for all $k \in \N$.
Example graphs $\gs$ associated with two time-varying delay models are shown in Figure~\ref{fig:logic}.
    % For example, Fig.~\ref{fig:logic} illustrates possible graphs $\gs$ associated with two time-varying delays models. 
    Figure~\ref{fig:logic:a} shows a graph that encodes delays  taking values $h_k\in \{2, \ldots, 5\}$ with a bounded rate of change  $|h_{k+1}-h_k|\leq 1$ for all $k \in \N$. The graph in Figure~\ref{fig:logic:b} describes the packet drop case as a reset system: the delay $h_k \in \{0, 1, 2, 3\}$ either increases by one every time step (packet drop) or drops to 0 (successful packet transmission). The maximum  possible delay is $3$ steps, upon which the next delay is 0.
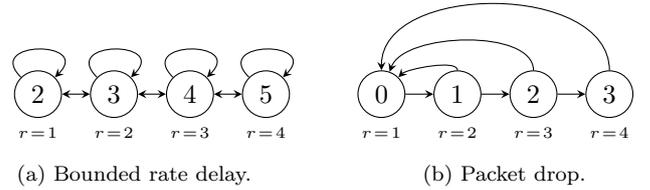
\begin{figure}[h]
    \centering
    \begin{subfigure}{0.45\linewidth}   
    \centering
        \begin{tikzpicture}[>=stealth, node distance=10mm, auto]
        % Define nodes
        \node[circle, draw] (0) {2};
        \node[circle, draw, right of=0] (1) {3};
        \node[circle, draw, right of=1] (2) {4};
        \node[circle, draw, right of=2] (3) {5};

        \node[below=0.0pt of 0] {\tiny $r\!=\!1$};
        \node[below=0.0pt of 1] {\tiny $r\!=\!2$};
        \node[below=0.0pt of 2] {\tiny $r\!=\!3$};
        \node[below=0.0pt of 3] {\tiny $r\!=\!4$};
        
        % Bidirectional edges between nodes
        \draw[<->] (0) -- (1);
        \draw[<->] (1) -- (2);
        \draw[<->] (2) -- (3);
        
        % Loops
        \draw[->, looseness=4, out=135, in=45] (0) to (0);
        \draw[->, looseness=4, out=135, in=45] (1) to (1);
        \draw[->, looseness=4, out=135, in=45] (2) to (2);
        \draw[->, looseness=4, out=135, in=45] (3) to (3);
        \end{tikzpicture}
        \caption{Bounded rate delay.}
        \label{fig:logic:a}
    \end{subfigure}
    \hfill
    \begin{subfigure}{0.45\linewidth}
    \centering
        \begin{tikzpicture}[>=stealth, node distance=10mm, auto]
        % Define nodes
        \node[circle, draw] (0) {0};
        \node[circle, draw, right of=0] (1) {1};
        \node[circle, draw, right of=1] (2) {2};
        \node[circle, draw, right of=2] (3) {3};

        \node[below=0pt of 0] {\tiny $r\!=\!1$};
        \node[below=0pt of 1] {\tiny $r\!=\!2$};
        \node[below=0pt of 2] {\tiny $r\!=\!3$};
        \node[below=0pt of 3] {\tiny $r\!=\!4$};
        
        % Bidirectional edges between nodes
        \draw[->] (0) -- (1);
        \draw[->] (1) -- (2);
        \draw[->] (2) -- (3);
        
        % Loops
        \draw[<-, looseness=1, out=90, in=90] (0) to (3);
        \draw[<-, looseness=0.75, out=70, in=90] (0) to (2);
        \draw[<-, looseness=0.5, out=45, in=90] (0) to (1);
        \end{tikzpicture}
        \caption{Packet drop.}
        \label{fig:logic:b}
    \end{subfigure}
    \caption{Examples of switching graphs $\mathcal{G}$ for different models of time-varying delay.}
    \label{fig:logic}
\end{figure}

If the delay can take the values of $h_k \in \{0, 1, 2\}$, the modes of this switched gradient descent scheme  may be described using the state space matrices
\begin{align}
    \mat{c|c}{\Acl_r & \Bcl_r \hl \Ccl_r & \Dcl_r} \in \left\{ \mat{ccc|c}{0 & 0 & I & I \\
    I & 0 & 0 & 0 \\ 
    0 & 0 & I & -\gamma I\hl
    0 & 0 & I & 0}, \ \mat{ccc|c}{0 & 0 & I & I \\
    I & 0 & 0 & 0 \\ 
    0 & 0 & I & -\gamma I\hl
    I & 0 & 0 & 0},  \ \mat{ccc|c}{0 & 0 & I & I \\
    I & 0 & 0 & 0 \\ 
    0 & 0 & I & -\gamma I\hl
    0 & I & 0 & 0} \right\}
 \label{eq:gd_example}
\end{align}
}

}

\subsection{Algorithms over Networks}

\rw{The setting of switched networked optimization algorithms includes an algorithm  \eqref{eq:alg_deployment} in which $\Sigma$ is formed through the interconnection of a network and a controller.}  For each fixed mode $r \in \{1, \ldots, N_s\}$, the network $P$ is described by 
\begin{subequations}\label{eq:optimization_algo:general}
\begin{align}
\label{eq:optimization_algo:general:plant}
    P_r : & & \mat{c}{
       x_{k+1}^{\rw{N}} \hl z_k \\ y_k
    }&= \mat{c|cc}{
        A_r & B_{r,1} & B_{r,2} \hl
        C_{r,1} & D_{r, 11} & D_{r,12} \\
        C_{r,2} & D_{r,21} & D_{r, 22}
    }\mat{c}{
      x_k^{\rw{N}} \hl w_k \\ u_k}.
    \intertext{The network $(P, \gs)$ will be interconnected with a switched controller $(K, \gs)$,  which is described by }
        K_r: & & \mat{c}{
            \xi_{k+1} \hl u_k
        } &= \mat{c|c}{
            A^K_r & B^K_r \hl
            C^K_r & D^K_r 
        }\mat{c}{
            \xi_k \hl y_k
        }\rw{.}\label{eq:optimization_algo:general:controller}
% \end{align}
\end{align}
\end{subequations}    
The resulting system is the interconnection $(\Sigma, \gs) = (P \star K, \gs)$.

Time-varying delays between the controller $K$ and the gradient $\nabla f$ can be modeled as a switched network $P$.
% can be modeled as switched systems
\rw{If the delay can take values between} $0$ and $h_{\max}$ \rw{inclusive, then the delay in a single coordinate} can be mapped to a switched system  \citep{conte2020modeling} with modes $r \in \{1, \ldots, h_{\max} \rw{+1\}}$ as
% where each delay $s_k$ is assigned one state-space realization
% \begin{subequations}

\begin{align}
\label{eq:switched_delay}
    \begin{aligned}
    {A}^d_{r} &=    \begin{pmatrix}
        0_{1 \times h_{\max}-1} & 0 \\
        I_{h_{\max}-1} & 0_{h_{\max}-1 \times 1}
    \end{pmatrix},  & 
    {B}^d_{r} &= \begin{pmatrix}
        1 \\ 0_{ h_{\max}-1 \times 1}
    \end{pmatrix}, \\
    {C}^d_{r} &= \begin{cases}
        0_{1 \times h_{\max}-1} & r = 1 \\
        e_{r-1}^\top & r \neq 1,
    \end{cases} & 
    {D}^d_{r} &= \begin{cases}
        1 & r = 1 \\
        0 & r \neq 1,
    \end{cases}
    \end{aligned}
\end{align}
in which the standard unit vector has dimensionality $e_{r-1} \in \R^{h_{\max}}$ for all $r \in \N$.
% \end{subequations} 
% The system in \eqref{eq:switched_delay} will flush out the initial condition $\xi^0$ after $h_{\max}$ time steps. 
A variable delay  \rw{taking} values in $h_k \in \{0,1,2\}$ \rw{over all coordinates }can be represented as the switched system
    \begin{align}
    \mat{c|c}{
            {A}^d_r & {B}^d_r \hl
            {C}^d_r & {D}^d_r
        } \in
      \left\{\left(\begin{array}{cc|c}     
        0 & 0 &  I \\
        I & 0 &  0 \\ \hline   
        0 & 0 &  I
        \end{array}
        \right), \left( \begin{array}{cc|c}     
        0 & 0 & I \\
        I & 0 & 0 \\ \hline 
        I & 0 & 0 
        \end{array}
        \right), \left( \begin{array}{cc|c}              
        0 & 0 & I \\
        I & 0 & 0  \\\hline
        0 & I  & 0
        \end{array}
        \right) \right\}. \label{eq:variable_delay_example}
    \end{align}

\rw{The gradient descent example in \eqref{eq:gd_example} can be obtained using the delay primitives in \eqref{eq:variable_delay_example} through the following  interconnection for each $r \in \{1, 2, 3\}$:
\begin{align}
\Sigma_r &= P_r \star K_r, & 
    \mat{c|c}{\Acl_r & \Bcl_r \hl \Ccl_r & \Dcl_r} &=  \mat{c|c:c}{
            {A}^d_r & 0 & {B}^d_r \hl
            {C}^d_r & 0 & {D}^d_r \hdl
            0 & I & 0             
        }  \star \mat{c|c}{I & -\gamma I \hl I & 0}.
\end{align}}
\section{Optimization Algorithm Analysis}
\label{sec:analysis}

We begin by \rw{providing a sufficient condition based on regulation theory under which switched optimization algorithms can converge to the optimal point $z^*$ \eqref{eq:problem}. We then use linear matrix inequality methods to upper-bound the exponential convergence rate of a given switched optimization algorithm.}

\subsection{Regulation of Switched Optimization Algorithms}

We first describe a regulation property \citep{francis1976internal, wen2008switched} that is sufficient to ensure that trajectories of \eqref{eq:alg_deployment} converge to an optimal point $z^*$.
% The class $\mathcal{S}_{m, L}$ contains the function $f(z) = \frac{m}{2}\norm{z}_2^2$. The interc

\begin{thm}
\label{thm:convergence}
    A sufficient condition for the interconnection \eqref{eq:alg_deployment} to converge to an optimal solution $z^*$ for any $s \in \text{Path}(\gs),$ $f \in \s_{m, L}$ 
    % \rw{under Assumption \ref{assum:paths} }
    is that there exist matrices $(\Theta_r)_{r=1}^{N_s}$ such that 
    \begin{subequations}
\label{eq:regulation_closed}  
            \begin{align}
 &\forall (r,r') \in \es:   & \Acl_r \Theta_{r} - \Theta_{r'}&= 0 , \label{eq:regulation_closed_invariant}  \\
 &\forall r \in \{1, \ldots,  N_{s}\}:  & \Ccl_r\Theta_r &= I_c, \label{eq:regulation_closed_kernel} 
    \end{align}    
    \end{subequations}
    and that the following \rw{switched} Lur\'e interconnection is asymptotically stable  for all $v_0 \in \R^n$,  $f^0 \in \s_{m, L}^0$, $s \in \text{Path}(\gs)$:
    \begin{subequations}
\label{eq:switch_exo_post_coord_f0}
\begin{align}
    v_{k+1} &= \Acl_{s_k} v_k + \Bcl_{s_k} \tilde{w}_k &     \tilde{w}_k &= \nabla f^0(\tilde{z}_k),  \\
    \tilde{z}_k &= \Ccl_{s_k} v_k.
\end{align}
% The algorithmic interconnection then satisfies $\lim_{k \rightarrow \infty} x_k
\end{subequations}
\end{thm}
% \begin{proof}
\textit{Proof:} \rw{In order to apply regulation theory, we use error coordinates to center the optimization problem about zero.
Given a fixed orbit $(\{x^*_r\}_{r=1}^{N_s}, z^*, w^*)$ that solves \eqref{eq:fixed_point} and a switching sequence $s \in \text{Path}(\gs)$, we define the error coordinates $(\tilde{x}, \tilde{z}, \tilde{w})$ as 
\begin{align}
    \tilde{x}_k& := x_k - x^*_{s_k}, & \tilde{z}_k&:= z_k - z^*, &  \tilde{w}_k &:= w_k - w^* = w_k.
\end{align}
The error $e$ and  error-objective $f^0$ are defined as
\begin{align}
    e_k &:= z_k - z^*  & f^0(\tilde{z}) &:=  f(\tilde{z} + z^*). \label{eq:error_coord}
\end{align}

The signal $e$ will be treated as the regulated output, while $\tilde{z}$ is the input to the nonlinearity $f^0$. 
% The equality $\tilde{z} = e$ is key to application of a linear regulation framework in the nonlinear optimization setting. 
The algorithmic interconnection in  \eqref{eq:alg_deployment} can then be expressed using a constant external disturbance $z^*$ as
% The system in \eqref{eq:alg_deployment } can be written in extended form as
\begin{align}
   \mat{c}{\tilde{x}_{k+1} \hdl z^* \hl \tilde{z}_k \\ e_k} &= \mat{c:c|c}{\Acl_{s_k}  & 0 & \Bcl_{s_k}\hdl
   0 & I & 0\hl
    \Ccl_{s_k} & -I & 0 \\
    \Ccl_{s_k} & -I & 0} \mat{c}{\tilde{x}_{k} \hdl z^* \hl \tilde{w}_k}, & \tilde{w}_k = \nabla f^0(\tilde{z}_k). \label{eq:alg_deployment_err}
\end{align}
If the interconnection in \eqref{eq:alg_deployment_err} satisfies $\lim_{k \rightarrow \infty}(\tilde{z}_k, \tilde{w}_k) = (0, 0)$ for any initial condition $\tilde{x}_0 \in \R^n$, disturbance  $z^* \in \R^c$, function $f^0 \in \mathcal{S}_{m, L}^0,$ and switching sequence $s \in \text{Path}(\gs)$, then the original algorithmic interconnection in \eqref{eq:alg_deployment} satisfies $\lim_{k \rightarrow \infty} z_k = z^*$ for any initial state $x_0 \in \R^n$ switching sequence $s \in \text{Path}(\gs)$, and function $f \in \mathcal{S}_{m, L}.$ 
% The system in \eqref{eq:alg_deployment_err} is a nonlinear interconnection
After introducing a switching-dependent new state $v_k := x_k - \Theta_{s_k} z^*$,
\eqref{eq:alg_deployment_err} becomes
\begin{align}
   \mat{c}{v_{k+1} \hdl z^* \hl \tilde{z}_k \\ e_k} &= \mat{c:c|c}{\Acl_{s_k} & \Acl_{s_k} \Theta_{s_k}  - \Theta_{s_{k+1}} & \Bcl_{s_k} \hdl
   0 & I & 0\hl
    \Ccl_{s_k} & \Ccl_{s_k} \Theta_{s_k}  - I& 0 \\
    \Ccl_{s_k} & \Ccl_{s_k} \Theta_{s_k} - I & 0} \mat{c}{v_{k} \hdl z^* \hl \tilde{w}_k}, & \tilde{w}_k = \nabla f^0(\tilde{z}_k).\label{eq:alg_deployment_errflip}
\end{align}
% through the use of the coordinate transformation
% \begin{align}
%     \mat{cc}{I & -\Theta_{s_{k+1}} \\ 0 & I} \mat{cc}{\Acl_{s_k} & 0 \\ 0 & I} \mat{cc}{I & \Theta_{s_{k}} \\ 0 & I} = \mat{cc}{\Acl_{s_k} & \Acl_{s_k} \Theta_{s_k} - \Theta_{s_{k+1}} \\ 0 & I}.
% \end{align}

If the equations in \eqref{eq:regulation_closed}   are satisfied, then \eqref{eq:alg_deployment_errflip} decouples  the exogenous $z^*$ input from the evolution of $(v, \tilde{z}, e)$ as 
\begin{align}
   \mat{c}{v_{k+1} \hdl z^* \hl \tilde{z}_k \\ e_k} &= \mat{c:c|c}{\Acl_{s_k} & 0 & \Bcl_{s_k}\hdl
   0 & I & 0\hl
    \Ccl_{s_k} & 0 & 0 \\
    \Ccl_{s_k} & 0 & 0} \mat{c}{v_{k} \hdl z^* \hl \tilde{w}_k}, & \tilde{w}_k = \nabla f^0(\tilde{z}_k).\label{eq:alg_deployment_errflip_2}
\end{align}

The robust stability result in \eqref{eq:switch_exo_post_coord_f0} ensures that $\lim_{k \rightarrow \infty} v_k = 0$,  and therefore, that  $\lim_{k \rightarrow \infty}  \tilde{z}_k = \lim_{k \rightarrow \infty}  e_k = \lim_{k \rightarrow \infty} \Ccl_{s_k} v_k = 0$ for all initial points $v_0 \in \R^n$, switching sequences $s \in \text{Path}(\gs),$ and functions $f^0 \in \sml$. This proves convergence of the switching optimization algorithm to $z^*$. \qed

}

\begin{remark}
\rw{Theorem~\ref{thm:convergence} can be embedded into the regulation problem for switched linear systems in the literature by fixing the function $f\in \s_{m,L}$ in \eqref{eq:alg_deployment} to the very special case $f(z)=\frac{1}{2}m\|z-z^*\|^2$.
%implying $f^0(z)=\frac{1}{2}m\|z\|^2$. 
The interconnection \eqref{eq:alg_deployment} with the error output $e_k:=z_k-z^*$ then 
reads as
\begin{equation}
x_{k+1}=(\Acl_{s_k}+m\Bcl_{s_k}\Ccl_{s_k}) x_{k}-m\Bcl_{s_k}z^*\te{and}e_k=\Ccl_{s_k}x_{k}-z^*.   \label{eq:ssysnom}
\end{equation}
Regulation of linear systems amounts to ensuring that \eqref{eq:ssysnom} is asymptotically stable for $z^*=0$ and satisfies
$\lim_{k\to\infty} e_k=0$ for all $z\in\R^c$. 

There exists a wide literature on guaranteeing uniform exponential
stability of switched systems based on common quadratic or switched Lyapunov functions \cite{Daafouz2002}, which are amenable to 
the convex design of switching dependent controllers \cite{de2012gain}. In particular, the work in \cite{de2012gain} forms the basis for synthesis in the current paper.

For LTI systems, the additional regulation property is assured through the existence of a solution of the related
regulator equation as known from the internal model principle \cite{francis1976internal,stoorvogel2000performance}.
The extension to switched systems is handled in \cite{zattoni2013output} and \cite{Zhao2020} by a 
common solution of the regulator equation; if applied for \eqref{eq:ssysnom}, this equation boils down to \eqref{eq:regulation_closed} with
$\Theta$ being independent of the node index $r$. Motivated by almost regulation for LPV systems \cite{korouglu2008lpv}, 
our result is based on the use of mode-dependent solutions of the regulator equation 
to reduce conservatism and to even guarantee robust regulation against structured nonlinear uncertainties 
$\nabla f$ with $f\in\s_{m,L}$. 

The use of path-dependent \cite{Lee2006} or path-complete \cite{Ahmadi_2014} Lyapunov functions as surveyed in \cite{Lee2020}
permits to arrive at necessary and sufficient conditions for uniform exponential stability (which are as well
amenable to controller synthesis) by lifting the switched system into a larger system based on $\ell$-length histories. To the best of our knowledge, an analogous result to ensure necessary and sufficient conditions for output regulation 
with path-dependent solutions of the regulator equation has not been considered in the literature.}
% \green{A comment here or later: A key contribution is to combine this with ZF multipliers and LPV synthesis to design algorithms?}
\end{remark}
%switched Lypaunov functions 
%\cite{Daafouz2002} %switched LF -> synthesis 
%\cite{de2012gain}  %LPV synthesis 

%Regulation fixed solutions
%\cite{zattoni2013output}
%\cite{Zhao2020}

%Regulation in LPV setting: Motivation
%\cite{korouglu2008lpv}

%\cite{Lee2020} %nice survey
%\cite{Ahmadi_2014} %path-complete Lyapunov functions 

%possible extension to multiple

%Dissipativity, LMIs and LPV:
%ECC95 paper of mine
%SIAM book contribution from 2000 of mine

\subsection{Exponential Convergence Rate Analysis}
Theorem \ref{thm:convergence} offers a condition for asymptotic \rw{convergence} to $z^*$. 
\rw{The speed of asymptotic convergence can be judged by its exponential convergence rate.
 The algorithmic interconnection in \eqref{eq:alg_deployment} is $\rho$-exponentially convergent for $\rho \in (0, 1)$ if \rw{there exists a $\gamma_z > 0$ such that} the following inequality is obeyed by any sequence $z$ of  generated by the  algorithm 
\begin{align}
    \norm{z_k - z^*}_2 \leq \gamma_z \rho^k   \norm{z_0 - z^*}_2 & & \forall k \in \N. \label{eq:exp_stability_z}
\end{align}

The system $(\Sigma, \gs)$ forms a $\rho$-exponentially convergent \textit{switched optimization algorithm} with respect to the function $f \in \mathcal{S}_{m, L}$ if  \eqref{eq:switched_opt_limit} and \eqref{eq:exp_stability_z} both hold. The system is $\rho$-exponentially convergent in the state $x$ to the fixed orbit $\{x^*_r\}_{r=1}^{N_s}$ if the following inequality is satisfied at any initial condition $x_0$ and switching signal $s \in \text{Path}(\gs)$ of \eqref{eq:alg_deployment}
\begin{align}
  \norm{x_k - x^*_{s_k}}_2 \leq   \rho^{k} \norm{x_0 - x^*_{s_k}}_2 & & \forall k \in \N. 
\end{align}
% if the sequence $\{z_k\}_{k \in \N}$ satisfies the exponential stability property in  \eqref{eq:exp_stability_z} for any switching sequence $s \in \text{Path}(\gs)$ and initial point $x_0 \in \R^n$. 
We note that 
$\rho$-exponential stability in the state $x$ to the fixed orbit $\{x^*\}_{r=1}^{N_s}$ implies $\rho$-exponential stability in $z$ to $z^*$, because  the following inequality holds for every $s \in \text{Path}(\gs)$ and $k \in \N$:
\begin{align}
    \norm{z_k - z^*}_2 &= \norm{\Ccl_{s_k} x_k - \Ccl_{s_k} x^*_{s_k}}_2  \nonumber\\
    & \leq \sigma_{\max}(\Ccl_{s_k})\norm{x_k - x^*_{s_k}}_2 \leq \sigma_{\max}(\Ccl_{s_k}) \rho^{k} \norm{x_0 - x^*_{s_k}}_2.\label{eq:exp_stability_xz}
\end{align}}
The exponential-convergence-rate analysis problem \rw{over $\s_{m, L}$ is as follows:}
\begin{prob}
    Given a function class $\s_{m, L}$, a switched system  $(\Sigma, \gs) $, and a rate $\rho>0$, certify if $(\Sigma, \gs)$  is a \rw{$\rho$-convergent switched optimization algorithm  for all} $f \in \mathcal{S}_{m, L}$.
    \label{prob:analysis}
\end{prob}

To solve Problem \ref{prob:analysis}, we  use a condition for $\rho$-exponential stability of \eqref{eq:alg_deployment}.
% based on the interconnection of an exponentially transformed system in \eqref{eq:switched_system} 
%
\begin{cor}
    A sufficient condition for the algorithm in \eqref{eq:alg_deployment} to be $\rho$-convergent to $z^*$ for all $f \in \mathcal{S}_{m, L}$ as in \eqref{eq:exp_stability_z} is if the condition \eqref{eq:regulation_closed} holds, and  the following system is bounded (Lyapunov stable):
%     The system in \eqref{eq:switched_system} is internally Lyapunov stable if there exists a $\gamma > 0$ such that 
% \begin{align}
%     \forall k \in \N: \norm{x_k}_2 \leq \gamma \norm{x_0}_2 \label{eq:lyap_stability}
% \end{align}
    % is Lyapunov stable in the sense of   for all $s \in \text{Path}(\gs), \ f^0 \in \s_{m, L}^0$:
    % in \eqref{eq:switch_exo_post_coord_f0} is $\rho$-stable
        \begin{subequations}
\label{eq:exp_stable}
\begin{align}
    \bar{v}_{k+1} &= (\rho^{-1} \Acl_{s_k}) \bar{v}_k + (\rho^{-1} \Bcl_{s_k}) \bar{w}_k,     & \bar{w}_k &= \rho^{-k} \nabla f^0(\rho^{k} \bar{z}_k) \label{eq:exp_stable_v} \\
    \bar{z}_k &= \Ccl_{s_k} \bar{v}_k.
\end{align}
\end{subequations}
\end{cor}
% \begin{proof}
This proof follows arguments from \citep{SCHERER2023robustexp} about exponential stability.
Let $\mathcal{T}_{\rho^{-1}}$ denote the $\rho$-exponential signal weighting map, transforming a sequence $v$ as 
\begin{align}
    \mathcal{T}_{\rho^{-1}}: (v_0, v_1, v_2, \ldots) \mapsto (v_0, \rho^{-1} v_1, \rho^{-2} v_2, \ldots).
\end{align}

Defining the exponentially  weighted signals $\bar{v} = \mathcal{T}_{\rho^{-1}} v, \bar{e} = \mathcal{T}_{\rho^{-1}}e, \bar{w} = \mathcal{T}_{\rho^{-1}} w,$ we note that boundedness of $\bar{v}$ as $\norm{\bar{v}_k}_2 \leq \gamma \norm{\bar{v}_0}_2$ over all $k \in \N, x_0 \in \R^n$ for some $\gamma>0$, implies exponential stability of the unweighted signal $v$  as $ \norm{{v}_k}_2 \leq  \gamma \rho^{-k} \norm{{v}_0}_2$  for all $k \in \N,  v_0 \in \R^n$ \citep{desoer1975feedback}. \qed

% and the $\rho$-weighted static but  time-dependent function
%     \begin{align}
%      \bar{F}(k, z)& := \rho^{-k} \nabla f^0(\rho^k z)
%      \end{align}
     
    % \urg{Fill this in. Arguments from \cite{SCHERER2023robustexp} on exponential stability, signal spaces}
% \end{proof}

\rw{We also use the bar notation to denote the $\rho$-exponential weighting of a switched linear system as in \eqref{eq:exp_stable}, such as $\bar{G}$ resulting from $G$ with the following state space descriptions at each $r \in 1, \ldots, N_s$:
\begin{align}
    G_r &= \mat{c|c}{\Acl_r & \Bcl_r \hl \Ccl_r & \Dcl_r} &  \bar{G}_r &= \mat{c|c}{\rho^{-1} \Acl_r & \rho^{-1} \Bcl_r \hl \Ccl_r & \Dcl_r}. \label{eq:exp_weighting_sys}
\end{align}
}

% The oracle $\nabla f \in \s_{m, L}$ satisfies a family of dissipation relations:
Next, we use families of dissipation inequalities satisfied by gradients of 
functions $f$ in  $\s_{m, L}^0$ 
\rw{to construct multiplier relaxations that can certify $\rho$-convergence of switched algorithms.}
%%%%%%%%%%%%%%%%%%%%%%%%%%%%%%%%%%%%%%%%%%%%%%%%%%%%%%%%%%%%%%%
% \subsection{Acronyms/Initialisms}
% \subsection{Optimization Oracles}
% \begin{thm[Theorem 5 of \citep{SCHERER2023robustexp}]
% \label{thm:pass_prop}
%      Functions in $f \in \s_{m, L}$ satisfy a property $\forall z_1, z_2 \in \R^d:$
% \begin{align}
%     \nabla f_m(z_1) ( \nabla f^L(z_1) - \nabla f^L(z_2)) \geq \mathcal{V}_f(z_1) - \mathcal{V}_f(z_2). \label{eq:pass_prop}
% \end{align}
% In particular, choosing $z_2 = z^*$ yields an interpretation of \eqref{eq:pass_prop} as a passivity property between $\nabla f_m, \nabla f^L$ for members of $\s_{m, L}$. The passivity property can be strengthened to generate a family of dynamic IQCs:
% Theorem \ref{thm:pass_prop} may be used to generate a family of valid dynamic IQCs:
% Functions in $\s_{m, L}^0$ satisfy a class of dynamic IQCs:
\begin{lem}[Lemma 5 of \citep{scherer2023optimization}]
\label{lem:exp_iqc}
    Consider a function $f^0 \in \s^0_{m, L}$, an exponential rate $\rho \in (0, 1],$ and a sequence of scalar coefficients $\{\lambda_\nu \}_{\nu=0}^{\nu_{\max}}$ satisfying
    \begin{align}
        \lambda_\nu \leq 0\quad \forall \nu \geq 1, \quad 
        \qquad \textstyle \sum_{\nu=0}^{\nu_{\max}} \rho^{-\nu} \lambda_\nu > 0 \label{eq:iqc_lambda}.
    \end{align}
    % \begin{subequations}        
    % \begin{subequations}    
    %  Define respective $m$ and $L$ transformations of the $\rho$-weighted nonlinearity $\bar{F}:$
    %  \label{eq:loop_transform}
    %  \begin{align}     
    %  \bar{F}_m(k, z)& := \bar{F}(k, z) - m z, \\ 
    %  \bar{F}^L(k,  z) &:= L z - \bar{F}(k, z).& 
    % \end{align}
    % \end{subequations}
    % The sequences 
    For any sequence $z$, define the exponentially weighted sequences $\bar{z}, \bar{p}, \bar{q}, \bar{g}$     as 
    \begin{subequations}
    \label{eq:exp_signals}
    \begin{align}
        \bar{z}_k &:= \rho^{-k} z_k, & \bar{q}_k &:= \rho^{-k} \nabla f^0(\rho^{k} \bar{z}_k) - m \bar{z}_k,\\
        \bar{p}_k &:= L \bar{z}_k - \rho^{-k} \nabla f^0(\rho^{k} \bar{z}_k), & \bar{g}_k &:= \textstyle \sum_{\nu=0}^k \lambda_\nu \bar{p}_{k-\nu}.
    \end{align}
    \end{subequations}
    Then the sequences $\bar{q}, \bar{g}$  satisfy the dissipation inequality for all $T \in \N, T>0$.
    \begin{align}
        \textstyle \sum_{k=0}^{T-1} \bar{q}_k^{\top} \bar{g}_k \geq 0. \label{eq:exp_passive}
    \end{align}
    % \end{subequations}
\end{lem}

The coefficients $\lambda$ in \eqref{eq:iqc_lambda} parameterize a \rw{Zames-Falb Finite Impulse Response} (FIR)  filter $\Psi(\lambda) = \sum_{\nu=0}^{\nu_{\max}} \lambda_\nu \mathbf{z}^{-\nu}$.
The choice of $\lambda_0 = 1$ and $\lambda_{\nu} = 0$ for all $\nu \geq 1$ leads to a passivity property (identity filter). 
Cases with $\lambda_{\nu} \neq 0$ for some $\nu \geq 1$ describe dynamic filters. 
% This paper will restrict analysis and synthesis to causal multipliers for simplicity of explanation, but we note that non-causal multipliers are admissible in this framework \cite{SCHERER2023robustexp}. % Generalizations of Lemma \ref{lem:exp_iqc} in the IIR multiplier case may be found in \urg{?}.

% For a given set of filter coefficients $\lambda$, 

A state space realization of the system $\Psi(\lambda)$ may be defined using fixed matrices $A_f,$ $B_f$,  and $\lambda$-affine maps $C_f\rw{(\lambda)}, D_f\rw{(\lambda)}$ as 
    \begin{align}
       \mat{c}{\psi_{k+1} \hl \bar{g}_k}
     = \mat{c|c}{
        A_f & B_f \hl
        C_f(\lambda) & D_f(\lambda)}
    \mat{c}{
        \psi_{k} \hl \bar{p}_k}. \label{eq:filter}
    \end{align}

% we denote $\Psi(\lambda): \bar{p} \rightarrow \bar{r}$ as the $\lambda$-affine mapping realizing the filter $\lambda$ as a linear system.

% In order to apply Lemma \ref{lem:exp_iqc} to prove exponential stability of the algorithm \eqref{eq:alg_deployment}, we must conduct \rw{a} signal transformation to match the form of  \eqref{eq:exp_signals}:

\rw{We conduct a signal transformation in order to render the conditions in Lemma \ref{lem:exp_iqc} in the form of a star product}.
\begin{align}
    \begin{pmatrix}
        \bar{p}_k \\ \bar{q}_k
    \end{pmatrix} &= \begin{pmatrix}
        LI & -I \\ I & -mI
    \end{pmatrix}\begin{pmatrix}
        \bar{w}_k \\ \bar{z}_k
    \end{pmatrix}, \ \begin{pmatrix}
        \bar{p}_k \\ \bar{w}_k
    \end{pmatrix} =\begin{pmatrix}
        (L-m)I & I \\ I & mI
    \end{pmatrix}\begin{pmatrix}
        \bar{q}_k \\ \bar{z}_k
    \end{pmatrix}. \label{eq:loop_transform_eval}
\end{align}

% The loop transformation of the static \eqref{eq:loop_transform_eval} by an upper Linear Fractional Transformation w.r.t. the exponentially weighted algorithm \eqref{eq:switched_system_alg} is well-defined \citep{szabo2018transformations}, forming the per-mode transformed system

\rw{Given a set of filter coefficitents $\lambda$ and a rate $\rho$, we form the system $\hat{G}^\lambda$ from a switched system $(\Sigma, \gs)$ as
% Analysis of the interconnection \eqref{eq:alg_deployment} at rate $\rho$ is based on verifying dissipation properties of the system $\hat{G}^\lambda$ with
\begin{align}
    \hat{\Sigma}^\lambda = \Psi(\lambda) \left[ \begin{pmatrix}
        (L-m)I & I \\ I & mI
    \end{pmatrix} \star \bar{G}\right], \quad \text{as described by } \mat{c|c}{
        \hat{\Acl}_r^\lambda & \hat{\Bcl}_r^\lambda \hl
        \hat{\Ccl}_r^\lambda & \hat{\Dcl}_r^\lambda
    }. \label{eq:j_conversion}
\end{align}
}

\begin{thm}
    % Under Assumptions \ref{assum:m_L}-\ref{assum:cycle}, 
    \rw{The} switched algorithm $(\rw{\Sigma}, \gs)$ from \eqref{eq:alg_deployment} is certified as  $\rho$-exponentially convergent for all $f \in \s_{m, L}$ \rw{if  the condition in \eqref{eq:regulation_closed} is satisfied, and there exists a filter length $\nu_{\max} \in \N$, filter coefficients $\lambda \in \R^{\nu_{\max}}$ satisfying \eqref{eq:iqc_lambda}, and matrices $M_r \in \psd^{n+\nu_{\max}}, \  M_r \succ 0 $ for each $r \in 1, \ldots, s$ such that the following LMIs are satisfied for all  $\forall (r, r') \in \es: $ }
    % for each $(r, r') \in \es:$

    \label{thm:alg_analysis}
% \begin{subequations}

 \begin{align}
 % \text{find}_{M, \lambda}  \quad & \forall r \in 1, \ldots, N_s: \ M_r \in \psd^{n+\nu_{\max}}, \  M_r \succ 0 \\ 
     & \quad [\star]^\top \begin{pmatrix}
         M_r & 0 \\ 0 &-M_{r'}
     \end{pmatrix}\begin{pmatrix}
         \hat{\Acl}_r^\lambda & \hat{\Bcl}_r^\lambda \\ I & 0
     \end{pmatrix} + [\star]^\top\begin{pmatrix}
         0 & I_c \\ I_c & 0
     \end{pmatrix} \begin{pmatrix}
         \hat{\Ccl}_r^\lambda & \hat{\Dcl}_r^\lambda \\
         0 & I_c
     \end{pmatrix} \prec 0  \label{eq:antipassivity_switched}          
 \end{align}
 % \end{subequations}
 \end{thm}
 % \begin{pf}
 \textit{Proof: }
A detailed proof follows by using arguments from \cite{scherer2023optimization, scherer2025tutorial}.
The equation in \eqref{eq:antipassivity_switched} \rw{can be interpreted as} strict antipassivity of the system in \eqref{eq:j_conversion}. Lemma \ref{lem:exp_iqc} details that \rw{gradients} of all functions in $\s^0_{m, L}$ are passive under a sequence of transformations via the dissipation relation \eqref{eq:exp_passive}. Because the interconnection between a strictly antipassive system and a passive system is bounded, feasibility of \eqref{eq:antipassivity_switched} ensures that the system in \eqref{eq:alg_deployment} is $\rho$-exponentially stable for all members of $\s^0_{m, L}$.  \qed
% Because the system in \eqref{eq:switched_system_alg} also satisfies the regulation equation \eqref{eq:regulation_closed}, 

% Theorem \ref{thm:convergence} ensures that stability under $\s_{m, L}^0$ interconnection, and the satisfaction of regulation \eqref{eq:regulation_closed} \rw{by the theorem statement} implies convergence to the fixed point $z^*$ for all $f \in \s_{m, L}.$ \qed
% The transformation in \eqref{eq:loop_transform_eval}

 % \urg{fix this}
 % The filter $\lambda$ is a valid dissipation relation for $\nabla f$ under Lemma \ref{lem:exp_iqc}. The transformation \eqref{eq:j_conversion} converts the system $(\Acl_r, \Bcl_r, \Ccl_r)$ into a form acceptable for passive interconnection with respect to \eqref{eq:exp_passive}. The inequality in \eqref{eq:antipassivity_switched} is an instance of antipassivity as reviewed in \eqref{eq:antipassive}. As such, feasibility of \eqref{eq:antipassivity_switched} certifies $\rho$-convergence of the algorithm.
     % This follows from standard arguments used in \urg{Does it? Fill this in}
 % \end{pf}

The problem in \eqref{eq:antipassivity_switched} is a finite-dimensional convex program for fixed $\rho$ and filter order $\nu_{\max}$. \rw{For fixed $\rho$, using an interior point method to solve the LMI system in \eqref{eq:antipassivity_switched} via semidefinite programming within a numerical tolerance of $\epsilon$ involves $O(\abs{\es}^{1/2} c (n+\nu_{\max}) \log(\epsilon^{-1}))$ steps,  with a per-step iteration complexity of $O(\abs{\es} c^6 (n+\nu_{\max})^6)$ \cite{ben2001lectures}. This worst-case complexity can be reduced by using solvers specialized for Kalman-Yakubovich-Popov-Lemma  problems \cite{vandenberghe2005interior}.}
% by fixing a maximal order $\nu_{\max}$ with $\forall \nu > \nu_{\max}: \ \lambda_\nu = 0$. 

Exponential convergence rates of $\Sigma$ can be upper-bounded by performing \rw{a} bisection on $\rho$. The accuracy of the bound can be increased by increasing $\nu_{\max},$ but there is no guarantee of convergence to the true rate bound. 

% \subsection{Computational Complexity}

\rw{
}

% Exponential performance of the algorithm $\tilde{G}$ can be guaranteed by performing bisection 

% The algorithm in \eqref{eq:switched_system_alg} is first 

% \subsection{Dynamic Multipliers}
% % Problem \ref{prob:analysis} will be approached by bisection o

% \urg{
% Discuss optimization algorithms. Define the class $S(m, L)$. Discuss IQCs associated with this class (starting with loop transformations). How this links to passivity and anti-passivity.

% Exponential convergence/rate bounds/exponential transformations.

% }
\section{Algorithm Synthesis}
% \urg{The main theory}
\label{sec:synthesis}

We now focus on the synthesis of switching optimization methods.
% The regulation property from Theorem \ref{thm:convergence} will be guaranteed through the introduction of an internal model. 

% The algorithm $(\tilde{A}_r, \tilde{B}_r, \tilde{C}_r)$ in \eqref{eq:switched_rystem_alg} could arise from the interconnection between a switched network model of 
% \begin{align}
%     P_r : & & \begin{pmatrix}
%         z \\ y
%     \end{pmatrix} &= \begin{pmat}[{|.}]
%         A_r & B_{r,1} & B_{r,2} \cr \-
%         C_{r,1} & D_{r, 11} & D_{r,12} \cr
%         C_{r,2} & D_{r,21} & D_{r, 22}\cr
%     \end{pmat} \begin{pmatrix}
%         w \\ u
%     \end{pmatrix} \\
%     \intertext{and a switched linear controller}
%         K^s: & & u &= \begin{pmat}[{|}]
%             A^c_r & B^c_r \cr \-
%             C^c_r & D^c_r \cr
%         \end{pmat}y.  \label{eq:controller}
% \end{align}
% The direct-feedthrough $D_{yu}$ connection is left as 0 for simplicity.
The algorithm synthesis problem with fixed filter $\Psi(\lambda)$ satisfying \eqref{eq:iqc_lambda} and convergence bound $\rho$ is as follows: 
\begin{prob}
    Given a function class $\s_{m, L}$, a network $(P, \gs) $ from \eqref{eq:optimization_algo:general:plant}, and \rw{filter coefficients} $\lambda \in \R^{\nu_{\max}}$ satisfying \eqref{eq:iqc_lambda}, \rw{synthesize} a controller \eqref{eq:optimization_algo:general:controller} $(K, \gs)$ such that $(P \star K, \gs)$ from \eqref{eq:alg_deployment} forms \rw{a} $\rho$-exponentially convergent algorithm when interconnected with any $f \in \mathcal{S}_{m, L}$.
\end{prob}

% We will use output regulation theory \citep{francis1976internal} in the framework switched systems \citep{liu2001output, conte2021disturbance} to create output-feedback controllers.

\subsection{Assumptions}

The regulation property from Theorem \ref{thm:convergence} will be guaranteed through the introduction of an internal model. This internal-model-based regulation technique is only possible if a regulator equation is solvable:
\rw{
\begin{lem}
If $G = P \star K$ is a convergent optimization algorithm for all $f \in \mathcal{S}_{m, L}$, and $\{\Theta_r\}_{r=1}^{N_s}$ are solutions to \eqref{eq:regulation_closed} with respect to the plant representations $\{P_r \star K_r\}_{r=1}^{N_s}$, then there exists a solution \rw{$\{\Pi_r, \Gamma_r\}_{r=1}^{N_s}$} of the following \textit{regulator equation}
% \begin{subequations}    
\begin{align}
\forall (r, r') \in \es: & & \mat{cc}{A_r & B_{r, 2} \\ C_{r, \rw{1}} & D_{r, 12}} \mat{c}{\Pi_{r}\\ \Gamma_r} &= \mat{c}{\Pi_{\rw{r'}} \\ I_c}.\label{eq:regulation_open}
% \begin{pmatrix}
%     A_r - I & B_{r, 2} \\
%     C_{r, 2} & D_{r, 12}
% \end{pmatrix}
% \begin{pmatrix}
%     \Pi \\ \Gamma
% \end{pmatrix}    &=  \begin{pmatrix}
%     0_{n \times 1} \\ -1
% \end{pmatrix}.
\end{align} 
% \end{subequations}
\label{lem:regulator_eq}
\end{lem}
\textit{Proof:}
     We can vertically partition each matrix  $\Theta_r$ as $\Theta_r = [\Pi_r^\top, \Xi_r^\top]^\top$ for every $r \in \{1, \ldots, N_s\}$. Given that the star product $P_r \star K_r$ requires that $(I -  D_2^K D_{r, 22})$ is invertible at each $r$, the following  matrices 
% \begin{subequations}

\begin{align}
\label{eq:exists_solution_reg}
    \Gamma_r &:= (I -  D_r^K D_{r, 22})^{-1}D^K_{r}(C_{2, r} \Xi_r +  C_{2, r} \Pi_r), \\
    \Phi_r &:= C_{2, r} \Pi_r + D_{2, r} \Gamma_r,
\end{align}
% \end{subequations}
exist and are unique at each $r \in \{1, \ldots, N_s\}$, and solve \eqref{eq:regulation_open}. \qed
}

 We \rw{therefore} require the following assumptions  to perform synthesis of switched $\rho$-convergent optimization algorithms that can be certified by Theorem \ref{thm:alg_analysis}:
\begin{assum}
\label{assum:regulation_open} Given a network $P$ in \eqref{eq:optimization_algo:general:plant}, the regulator equation \eqref{eq:regulation_open} has a solution $\{\Pi_r, \Gamma_r\}_{r=1}^{N_s}.$
\end{assum}

\begin{assum}
    The matrix  $(I - D_{r, 11} m)$ is invertible  for each $r \in \{1, \ldots, N_s\}$. \label{assum:wellposed}
\end{assum}

% & & 
% A_{r} \Pi_r + B_{r, 2} \Gamma_r &= \Pi_{r'} \\
% & & C_{r, \rw{1}} \Pi_r + D_{r, 12} \Gamma_r &= -\rw{I_c}.

% of well-posedness of a feedback interconnection:

% Validity of Assumption \ref{assum:regulation_open} allows us to construct a controller such that the regulation property \eqref{eq:regulation_open} will be satisfied by construction. 

\rw{If Assumption \ref{assum:regulation_open} is violated, then there cannot exist a controller $K$ such that  $P \star K $ can  be certified as convergent by Theorem \ref{thm:convergence} (Lemma \ref{lem:regulator_eq}). Solvability of \ref{assum:regulation_open} is a property only of the network $P$, and is indepenedent of $K$.}
% \rw{Assumption \ref{assum:regulation_open} is also a necessary condition for a switched optimization algorithm to . 
% }
\rw{Assumption \ref{assum:wellposed} ensures the interconnection of $D_{r, 11}$ and the error system with the known function $f^0(\tilde{z}) = \frac{m}{2}\norm{\tilde{z}}_2^2, \ \nabla f^0(\tilde{z}) = m \tilde{z}$ as described by $(mI) \star P_r$ is well-posed for each mode $r$. This well-posedness assumption may be connected to the linear regulation condition \eqref{eq:ssysnom}.}

\subsection{Regulation}

% In order to apply the regulation framework to generate controllers, we make an assumption based on the solvability of an open-loop regulation equation:

% \begin{figure}[h]
%     \centering
% \input{img/interconnections/ic_internal_model_control}
%     \caption{Internal model control.}
%     \label{fig:ic_internal_model_control}
% \end{figure}
Given a solution $\{\Pi_r, \Gamma_r\}_{r=1}^{N_s}$ to \eqref{eq:regulation_open}, we define  $\Phi_r := C_{r, 2} \Pi_r + D_{r, 22} \Gamma_r$ \rw{for each $r \in \{1, \ldots, N_s\}.$} \rw{We construct controllers $K$ as the interconnection between fixed internal models $Q$ and designed subcontrollers $R$.}
% \urg{Need to get the proofs for this}
An internal model $Q$ inspired by \citep{stoorvogel2000performance} (originally derived for a non-switched system) is described for each $r \in \{1, \ldots, N_s\}$ by
\begin{align}
    Q_r: & & 
    \mat{c}{
        \omega_{k+1} \hl u_k \\ \tilde{y}_k}&= \mat{c|c:cc}{
        \rw{I_c} & 0 & I & 0 \hl
        -\Gamma_r & 0 & 0 & I  \\
        \Phi_r & I & 0 & 0}
     \mat{c}{        
       \omega_{k} \hl y_k \hdl \tilde{u}_{1k} \\ \tilde{u}_{2k}}. \label{eq:internal_model}
\end{align}
\rw{The subcontrollers $R$ are described by }
% We restrict the control synthesis task to the design of subcontrollers $R_r$ in interconnection with $P_r \star Q_r$:
\begin{align}
    R_r: & & 
        \mat{c}{
            \tilde{\xi}_{k+1} \hl
            \tilde{y}_k
        }
     &= \mat{c|c}{
        A_{r, c} & B_{r, c} \hl
        C_{r, c1} & D_{r, c1} \\
        C_{r, c2} & D_{r, c2} 
    } \mat{c}{
        \tilde{\xi}_{k} \hl \tilde{u}_{1, k} \\ \tilde{u}_{2, k}
    }, \label{eq:subcontroller}
\end{align}
\rw{for which the controller $K$ in  \eqref{eq:optimization_algo:general:controller} can be expressed as the interconnection $K = Q \star R$. This interconnection may be described for each $r$ as}
\begin{align}
\label{eq:control_association}
    K_r: \ \mat{c|c}{
            A^K_r & B^K_r \hl
            C^K_r & D^K_r 
        } = \mat{c|c:cc}{
            I + D_{r, c1} \Phi_r & C_{r, c1} & D_{r, c1} \hl
            B_{r, c} \Phi_r & A_{r, c} & B_{r, c} \\
            -\Gamma_r + D_{r, c2} \Phi_r & C_{r, c2} & D_{r, c2} }. 
\end{align}

% \input{img/interconnections/ic_internal_model_control}

% In order to ensure that the ora
% The following structur\rw{al} constraint imposed on the controller $R$ at each $r \in \{1, \ldots, N_s\}$ ensures that there are no algebraic loops in the evaluation of $\nabla f$:

\rw{Requiring that $D_{r, 11} +  D_{r, 12} D_{r, c2} D_{r, 21} = 0$ at all modes $r \in, 1\ldots, s$ when designing the subcontroller $R$ ensures that $P \star Q \star R$ has  no algebraic loops in the evaluation of $\nabla f$.}
% \begin{align}
%     \label{eq:no_dependence}
% \end{align}
% % is added to 

% that the function $f$ is only observed through the oracle $\nabla f$, with no algebraic loops or proximal operator evaluations in the algorithm.

\begin{thm}
\label{thm:regulation_control}
    If Assumption \ref{assum:regulation_open} is satisfied, then 
    the controller structure \eqref{eq:control_association} ensures that $P \star Q \star R$ satisfies condition \eqref{eq:regulation_closed}
    for output regulation from Theorem \ref{thm:convergence}.
    \end{thm}
% \begin{pf}
\textit{Proof:}
% We focus on the case of $D_{r, 22} = 0$ (no direct feedthrough) to simplify the proof.
This follows from standard arguments in regulator theory from \citep{francis1976internal} with respect to the model structure in \citep{stoorvogel2000performance}.
Describing $P_r \star Q_r \star R_r$ as
\begin{align}
\mathcal{A}_{r} &= \mat{cc|c}{
        \mathcal{A}_r + B_{r, 2} D_{r, c2} C_{r, 2}& -B_{r, 2} (\Gamma_r -D_{r, c2} \Phi_r) & B_{r, 2} C_{r, c2} \\ 
    D_{r, c1} C_{r, 2} & I + D_{r, c1} \Phi_r & C_{r, c1} \hl
        B_{r, c} C_{r, 2} & -B_{r, c} \Phi_r  & I-A_{r, c}
}\nonumber \\
\mathcal{B}_{r} &=   \begin{pmatrix}
B_{r, 1}\label{eq:internal_model_network_connected} \\
   D_{r, c1} D_{r, 11} \\
       B_{r, c} D_{r, 11}\end{pmatrix}, \ 
       \mathcal{C}_{r} = \begin{pmatrix}
           (C_{r, 1} + D_{r, 12} D_{r, c2} C_{r, 2})^\top \\ (-D_{r, 12} (\Gamma_r - D_{r, c2} \Phi_r))^\top \\ (D_{r, 12} C_{r, c2})^\top
       \end{pmatrix}^\top, \   \mathcal{C}_{r} = 0,
%      \begin{pmat}[{, \ldots, |}]
%     A_r + B_{r, 2} D_{r, c2} C_{r, 2}& -B_{r, 2} (\Gamma -D_{r, c2} \Phi_r) & B_{r, 2} C_{r, c2} & B_{r, 1} \cr 
%     D_{r, c1} C_{r, 2} & I + D_{r, c1} \Phi_r & C_{r, c1} & D_{r, c1} D_{r, 11} \cr
%         B_{r, c} C_{r, 2} & -B_{r, c} \Phi_r  & I-A_{r, c} &  B_{r, c} D_{r, 11} \cr\-
%     C_{r, 1} + D_{r, 12} D_{r, c2} C_{r, 2} & -D_{r, 12} (\Gamma - D_{r, c2} \Phi_r)& D_{r, 12} C_{r, c2} &  0 \cr
% \end{pmat}.
% \end{align}
% The relevant matrix in Assumption \ref{assum:regulation_open} is
% \begin{align}
%     \begin{pmatrix}
%         I - A_r \\
%         -C_r
%     \end{pmatrix} &= \begin{pmatrix}
%             I - A_r - B_{r, 2} D_{r, c2} C_{r, 2} & B_{r, 2} (\Gamma - D_{r, c2}\Phi) & -B_{r, 2} C_{r, c2}  \\
%     -D_{r, c1} C_{r, 2} & -D_{r, c1} \Phi_r & -C_{r, c1} \\
%         -B_{r, c} C_{r, 2} & B_{r, c} \Phi_r  & -A_{r, c} \\
%     -C_{r, 1} - D_{r, 12} D_{r, c2} C_{r, 2}& D_{r, 12} (\Gamma - D_{r, c2}\Phi) & -D_{r, 12} C_{r, c2} \cr
%     \end{pmatrix}
% % \end{align}
\end{align}
\begin{align}
\intertext{and choosing  $\Theta_r = [-\Pi_r^\top, I, 0]^\top$ and $\Theta_{r'} = [-\Pi_{r'}^\top, I, 0]^\top$ yields }
% \begin{align}
     \begin{pmatrix}
         \mathcal{A}_r \\
        \mathcal{C}_r
    \end{pmatrix} \Theta_r  - \Theta_{r'}&= \mat{c}{
        \Pi_{r'} - A_r \Pi - B_{r, 2} \Gamma_r \\
        -D_{r, c1} (C_{r, 2}\Pi_r - ( C_{r, 2} \Pi_r)) \\
        -B_{r, c1} (C_{r, 2}\Pi_r - (C_{r, 2} \Pi_r)) \hdl
        -(C_{\rw{r, 1}} \Pi_r + D_{r, 12} \Gamma_r)
    } = \mat{c}{0 \\ 0 \\ 0 \hdl I_c}, \label{eq:regulation_many}
    \end{align}
which certifies \eqref{eq:regulation_closed}.
% The final step \rw{in} \eqref{eq:regulation_many} is taken by  applying relation \eqref{eq:regulation_open} and constraint $D_{r, 11} +  D_{r, 12} D_{r, c2} D_{r, 21} = 0$ for each pair $(r, r') \in \es$. 
\qed
% which proves the Lemma.
% \end{pf}

% \begin{rem}
%     The per-mode integrator-involved internal model in \eqref{eq:internal_model} is a sufficient but not necessary condition for regulation by Lemma \ref{lem:regulation_control}. As noted in Remark \ref{rmk:core_presence}, per-mode integrators are only necessary in modes $r$ with self-loops.
% \end{rem}
% \urg{Something about how the controller in \eqref{eq:control_association} achieves regulation, and satisfies }

% \begin{exmp}
%     Systems whose network dynamics are characterized by per-channel time-varying delay  with maximal $H_1^{\max}$ before $\nabla f$ and $H_2^{\max}$ after $\nabla f$ will satisfy Assumption \ref{assum:regulation_open} with an $r$-independent solution of 
%     \begin{align}
%         \Pi &= \begin{pmatrix}
%             1_{1 \times H_1^{\max}} & 0_{1 \times H_2^{\max}}
%         \end{pmatrix} \otimes \rw{I_c}& \Gamma &= -\rw{I_c} & \Phi = 0.
%     \end{align}    
% The resultant $r$-independent internal model is an  integrator which reads
% \begin{align}
%    Q_r:  \mat{c}{
%         \omega_{k+1} \hl u_k \\ \tilde{y}_k}
%     =\mat{c|c:cc}{
%         \rw{I_c} & 0 & \rw{I_c} & 0 \hl
%         \rw{I_c} & 0 & 0 & \rw{I_c} \\
%         0 & \rw{I_c} & 0 & 0 
%     }  \mat{c}{
%         \omega_k \hl y_k \hdl \tilde{u}_{1k} \\ \tilde{u}_{2k}}.
% \end{align}

% This model involves a proper integrator $(\tilde{u}_1, \tilde{u}_2) \rightarrow u$ and a direct output connection $\tilde{y} \rightarrow y$.
    
% \end{exmp}

\subsection{Controller Synthesis}

% The controller synthesis problem can be posed as finding a sub-controller $R_r$  ensuring the interconnection between $R_r$ and the following generalized plant is stable when connected to $\nabla f$:
% in a generalized plant framework, as visualized in Figure \ref{fig:ic_generalized_plant} with a per-mode central plant 

% The interconnection of network dynamics \eqref{eq:optimization_algo:general:plant} and the internal model \eqref{eq:internal_model} for each $r \in \{1, \ldots, N_s\}$ is given by $G_r = P_r \star Q_r$, as described by 
% \begin{align}
% \label{eq:plant_model_conn}
%     G_r: & \mat{c}{
%         \tilde{x}_{k+1} \\ \omega_{k+1} \hl z_k \\ \tilde{y}_{k}
%     } &= \mat{cc|c:cc}{
%         A_{r} & B_{r, 1} & B_{r, 1} & B_{r, 2} & \Pi_{r}\\ 
%         0 & \rw{I_c} & 0 & 0 & \rw{I_c} \hl
%         C_{r, 1} & D_{r, 11} & D_{r, 11} & D_{r, 12} & 0 \\
%         C_{r, 21} & D_{r, 21} & D_{r, 21} & D_{r, 22} & 0
%     }\mat{c}{
%         \tilde{x}_{k} \\ \omega_k \hl w_k \hdl \tilde{u}_{1, k} \\ \tilde{u}_{2, k}
%     }.
% \end{align}
% The system in \eqref{eq:plant_model_conn} has a description of
% % The connection of the network dynamics and the per-mode internal model is
% \begin{align}
%     G_r = P_r \star Q_r = \begin{pmat}[{|}]
%         \hat{A}_r & \hat{B}_r \cr \- \hat{C}_r & \hat{D}_r \cr
%     \end{pmat}.
% \end{align}
The controller synthesis problem can be posed as finding a sub-controller $(R, \gs)$  ensuring that $(P \star Q \star R, \gs)$ is a $\rho$-convergent optimization algorithm for any $f \in \mathcal{S}_{m, L}$. \rw{For a fixed filter $\lambda$ satisfying \eqref{eq:iqc_lambda}, the control synthesis problem may be solved by finding per-mode controllers $\{R_r\}_{r=1}^{N_s}$ such that the system $\hat{G}^\lambda$ with individual modes defined for each $r \in \{1, \ldots, N_s\}$ as
\begin{align}
    \hat{G}^\lambda_r := \Psi(\lambda) \left[ \begin{pmatrix}
        (L-m)I & I \\ I & mI
    \end{pmatrix} \star \bar{P}_r \star \bar{Q}_r \star \bar{R}_r\right],
\end{align}
satisfies the antipassivity LMI in \eqref{eq:antipassivity_switched} across all arcs $(r, r') \in \es$.
}
% \begin{figure}[h]
%     \centering
%     \input{img/interconnections/ic_generalized_plant}
%     \caption{Generalized plant.}
%     \label{fig:ic_generalized_plant}
% \end{figure}
% Control synthesis therefore requires the creation of a controller
% The IQC processing of the plant $G_r$ per Theorem \ref{thm:pass_prop} is 
% The anti-passive formulation for controller synthesis given $G_r, \rho, \lambda$ is to ensure that $R_r$ from \eqref{eq:subcontroller} renders 
% a $\lambda$-filtered, loop-transformed, and exponentially weighted plant anti-passive (detailed in Appendix A).
% the following plant anti-passive \eqref{eq:antipassive} with respect to switching sequences consistent with $\gs$:
% \begin{align}
%     \bar{G}^\lambda_r &= \Psi(\lambda) \cdot \left( \begin{pmatrix}
%         L & -1 \\ 1 & m
%     \end{pmatrix} \star \begin{pmat}[{|}]
%         \rho^{-1} \hat{A}_r & \rho^{-1} \hat{B}_r \cr \- \hat{C}_r & \hat{D}_r \cr
%     \end{pmat} \right).  \label{eq:plant_lambda}
% \end{align}
% satisfies the antipassivity relation  
% This antipassivity-restricted synthesis will be accomplished through the creation of $r$-dependent storage functions and nonlinear controller transformations. 
The controller synthesis and reconstruction procedure follows established methods for  Linear Parameter Varying synthesis, see  \citep{de2012gain} for further detail. The alternating search between $\lambda$ and $R$ for controller synthesis is described in Algorithm \ref{alg:algorithm1}. \rw{At the end of IterMax iterations, Algorithm \ref{alg:algorithm1} will return an exponential convergence rate $\rho$ and a subcontroller $R$ such that $(P \star Q \star R, \gs)$ is a $\rho$-convergent switched optimization algorithm.}

% , bisection tolerance $\rho_{\text{tol}}$. a $\rho$ range $[\rho_{\min}, \rho_{\max}]$

\begin{algorithm}[h]
% \begin{figure}
    % \centering

    \caption{Synthesis of Switched Optimization Algorithms}\label{alg:algorithm1}
    \begin{algorithmic}[1]
    \Require A switching graph $\gs$, a networked system $P$  \eqref{eq:optimization_algo:general:plant}, constants $(m, L)$, filter length $\nu_{\max}$, number of iterations IterMax.
    \State iter $\leftarrow$ 1
    \State Initial identity filter $\Psi(\lambda) \leftarrow 1$
    \State Compute internal models $Q$ from \eqref{eq:internal_model}
    \For{iter $\in 1, \ldots, $ IterMax }
    \State \textbf{Synthesis: } $R \leftarrow$ feasible subcontroller  with infimal $\rho$ under fixed filter $\Psi(\lambda) $ given $P \star Q.$
    % Apply bisection on $\rho$ to Theorem \ref{thm:alg_synth} with fixed multiplier $\lambda$ to obtain controller $P_r \star Q_r \star R_r$
    \State Recover the switched algorithm $(P \star Q \star R)$
    \State \textbf{Analysis:} $\lambda \leftarrow $ filter coefficients from  Theorem \ref{thm:alg_analysis} for analysis of  $(P \star Q \star R)$ with infimal $\rho$
    \EndFor \\
    \Return rate $\rho$, filter $\lambda$, subcontroller $R$
\end{algorithmic}
    % \includegraphics[width=0.5\linewidth]{}
    % \caption{Caption}
    % \label{fig:placeholder}
% \end{figure}
\end{algorithm}

\section{Numerical Examples}

\label{sec:examples}

Numerical experiments are implemented in MATLAB (2024a). All LMIs are solved using LMIlab \citep{gahinet1993lmilab} from the Robust Control Toolbox of MATLAB. The code to generate experiments is publicly available\footnote{\texttt{https://github.com/jarmill/time\_var\_opt}}. These examples aim to find the optimal value of the following 
 quadratic plus log-sum-exp function parameterized by terms $\Lambda \in \R^{c \times c}$, $b \in \R^c$, and $\ell \in \R$.
\begin{align}
    f(z; \Lambda, b, \ell) = &\frac{1}{2} z ^\top \Lambda z + b^\top z \label{eq:function_class} +  \ell \textstyle \log\left(\sum_{i=1}^d \rw{e}^{z_i}\right), 
\end{align}
\rw{The optimal value of \eqref{eq:function_class} must be found when the oracle $\nabla f(\cdot; \Lambda, b, \rw{\ell})$ is only accessible \rw{through} a time-varying network.}
Scalars $0 < m< L < \infty$ are given to define the function class \eqref{eq:function_class}.
The vector  $b$ is randomly drawn in $\R^c$. The matrix $\Lambda$ is a randomly generated symmetric matrix with eigenvalues between $m$ and $L'$ for a given $L' \in (m, L)$. The scalar $\ell$ is defined as $\ell:= L - L'$. Under these choices of parameters $(b, \Lambda, \ell)$, every instance of $f$ from \eqref{eq:function_class} is a member of $\mathcal{S}_{m, L}$. 
% The function class \eqref{eq:function_class} is contained in .

\subsection{Inhomogenous Networked System}

\begin{figure}[h]
    \centering
        % \begin{subfigure}{0.45\linewidth}
        \begin{tikzpicture}[>=stealth, node distance=10mm, auto]
        % Define nodes
        \node[circle, draw] (1) {1};
        \node[circle, draw, right of=1] (2) {2};
        \node[circle, draw, below=1] (4) {4};
        \node[circle, draw, right of=4] (3) {3};
        
        % Bidirectional edges between nodes
        \draw[->] (4) -- (1);
        \draw[->] (1) -- (2);
        \draw[->] (2) -- (3);
        \draw[->] (3) -- (4);
        
        % Loops
        \draw[->, looseness=4, out=135, in=225] (4) to (4);
        \draw[->, looseness=4, out=135, in=225] (1) to (1);
        \draw[->, looseness=4, out=45, in=315] (2) to (2);
        \draw[->, looseness=4, out=45, in=315] (3) to (3);
        \end{tikzpicture}
        % \caption{Bounded rate delay.}
        
    % \end{subfigure}
    \caption{Ring-structured switching graph} \label{fig:logic:ring}
    % \label{fig:placeholder}
\end{figure}
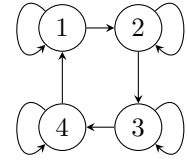

We first consider minimization of $f$ \eqref{eq:function_class} over a network governed by $N_s=4$ switching modes. The switching logic is displayed in Figure  \ref{fig:logic:ring}.

% \urg{Make a figure of this. Self loops, and cycle $1 \rightarrow 2 \rightarrow 3 \rightarrow 4$}.

\rw{The modes $\{P_r\}^4_{r=1}$ of the network $P$ are described by 
% The \rw{mapping per-mode systems   $(w, u) \rightarrow (z, y)$ are state-space realizations of the following transfer function matrices}
\small{
\begin{equation}
\label{eq:ring_network}
\begin{aligned}
    P_1 &: \mat{cc|cc}{0 & 0 & 1 & 0\\0 & 0.2 & 0 & 1 \hl 0 & 1 & 0 & 0 \\ -1 & 0 & 0 & 0}\otimes I_c, &   P_2 &: \mat{cc|cc}{0.2 & 0 & 0.25 & 0\\0 & 0.9 & 0 & 1 \hl 0 & 1 & 0 & 0 \\ -0.4 & 0 & -0.5 & 3}\otimes I_c,  \\
    P_3 &: \mat{cc|cc}{-0.3 & 0 & 0.5 & 0\\0 & -0.5 & 0 & 1 \hl 0 & 1 & 0 & 0 \\ -0.3 & 0 & 0.5 & 1}\otimes I_c,  &   P_4 &: \mat{cc|cc}{1.2 & 0 & 1 & 0\\0 & -0.2 & 0 & 1 \hl 0 & 1 & 0 & 0 \\ 1 & 0 & 0 & 2}\otimes I_c.
\end{aligned}
% \begin{aligned}
%     P_1 &= \begin{pmatrix}
%         0 & \frac{1}{\bz-0.2} I_c \\-\frac{1}{\bz} I_c & 0
%     \end{pmatrix} &   P_2 &= \begin{pmatrix}
%         0 & \frac{1}{\bz-0.9} I_c \\-\frac{1}{2(\bz-0.2)} I_c & 3 I_c
%     \end{pmatrix} \\
%     P_3 &= \begin{pmatrix}
%         0 & \frac{1}{\bz+0.5} I_c \\-\frac{1}{2(\bz-0.3)} I_c & 0
%     \end{pmatrix} &   P_4 &= \begin{pmatrix}
%         0 & \frac{1}{\bz+0.2} I_c \\-\frac{1}{\bz-1.2} I_c & 3 I_c
%     \end{pmatrix},
% \end{aligned}
\end{equation}}
}
\normalsize
% for which the full network dynamics are described by $P_r = P_r^1 \otimes I_c \ \forall r \in 1..N_s.$
% Applying gradient descent $(\alpha = 2/(L+m))$ is certified at $\nu_{\max}=4$ to have a convergence rate $\rho < 1.5974$. Triple momentum \citep{van2017fastest} with parameters tuned for the 0-delay setting has a convergence rate of $\rho < 1.6025$. In this experiment however, both algorithms fail to converge. 
The dynamics in \eqref{eq:ring_network} are chosen such that unstable channel dynamics are present in the subsystems (in $P_4$), and 
 that \eqref{eq:regulation_open} admits a solution with non-identical $\Pi$ and $\Gamma$ values between modes, \rw{with the specific values of}
 % and to demonstrate stabilization under unstable channel dynamics.
% Mode $r=4$ has a an unstable pole of $\bz = 1.2$ in its network description. 
% The internal model \eqref{eq:internal_model} constructed from the solution to the regulator equation is parameterized by 
\begin{align}
    \begin{pmatrix}
        \Gamma_r \\ \Phi_r
    \end{pmatrix}_{r=1}^4 = \begin{pmatrix}
        -0.8 & -0.1 & -1.5 & -1.2 \\
        0 & -0.3 & -1.5 & -2.4
    \end{pmatrix} \otimes I_c, \label{eq:regulator_1_solution}
\end{align}
and $\Pi_r = [0,  1]^\top \otimes I_c$ for all $r \in 1..4$. The solutions $(\Pi, \Gamma, \Phi)$ in \eqref{eq:regulator_1_solution} are unique. \rw{Thus, the method from \cite{zattoni2013output} requiring mode-independent solutions of the regulator equations cannot be used to certify regulation.}

% . Requiring $\Gamma_r = \Gamma^0$ at each $r$ for some constant $\Gamma^0$ fails Assumption \ref{assum:regulation_open}. 

\begin{figure}[!h]
    \centering
    \includegraphics[width=\linewidth]{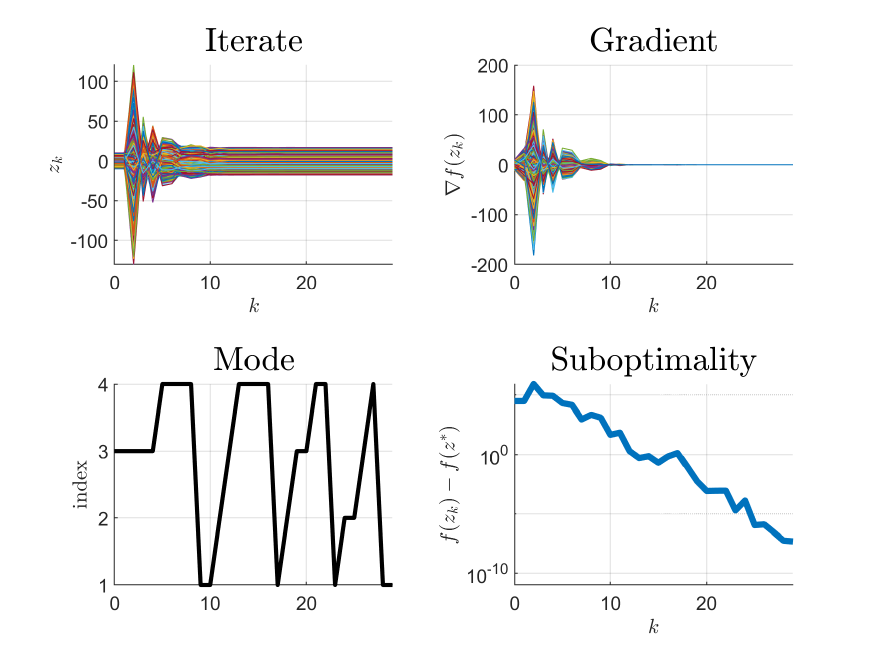}
    \caption{Algorithm execution with the $\rho =0.8913$ controller}
    \label{fig:lse_ring}
\end{figure}

We run Algorithm \ref{alg:algorithm1} for 3 iterations with Zames-Falb multipliers of order $\nu_{\max} = 3$. With only identity multipliers ($\lambda_0 = 1, \lambda_{\nu} = 0 \ \forall \nu \geq 1$), a controller $R$ of order 3 in each mode $r$ is obtained yielding $\rho \leq 0.8944$. At the end of 3 iterations, the analysis routine certifies a  controller $R$ of order 6 in each mode $r$ with $\rho \leq 0.8913$ under the multiplier $(\lambda_0, \lambda_1, \lambda_2, \lambda_3) = (1, -0.3922, -0.2977, -1.403\times 10^{-6})$.
% \end{align}
Figure \ref{fig:lse_ring} plots \rw{traces of optimizing}  an instance of $f$ over the network in \eqref{eq:ring_network} with $L'= 1.25, L= 1.6, \rw{c=400}$ by the developed controller. \rw{When starting at a randomly drawn initial state, the iterate $z$ converges to a nonzero value within 30 time steps (top left), and the gradient converges to zero (top right). The switching signal $s_k$ (bottom-left) arises from traversing the graph in \ref{fig:logic:ring} in a random walk choosing an edge uniformly at random. The suboptimality in the function value (bottom-right) shows that $f(z_{30}) - f(z^*) \approx 10^{-7}$ after 30 time steps.}

\begin{figure}[h]
    \centering
    \includegraphics[width=\linewidth]{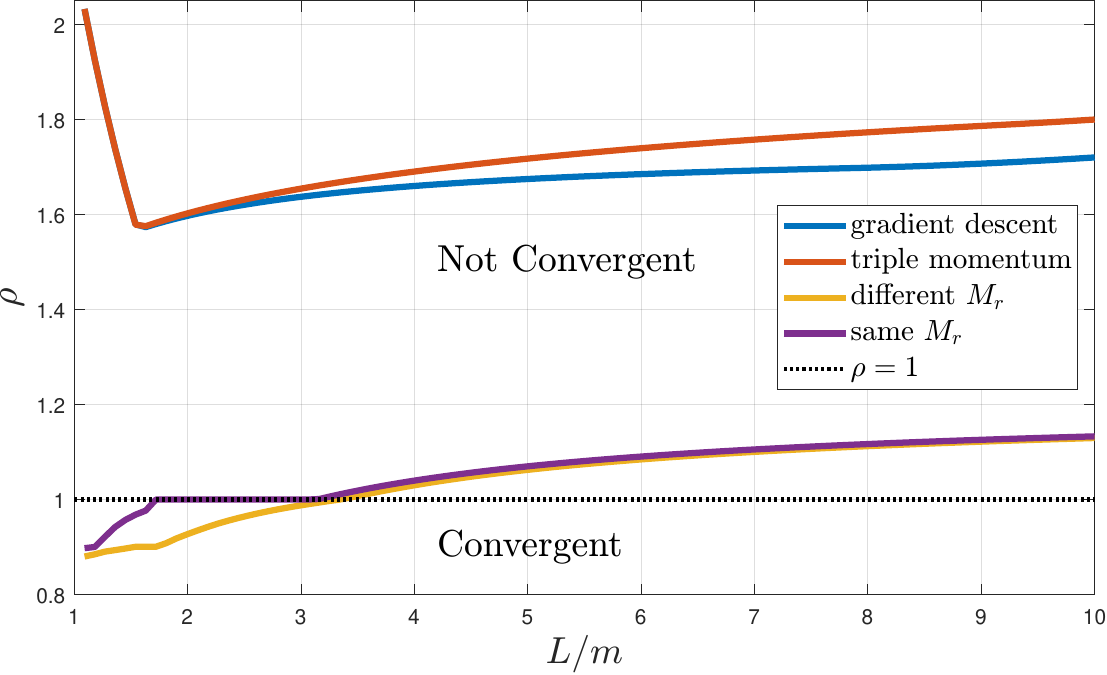}
    \caption{Analysis $\nu_{\max} = 3$ upper bounds on $\rho$}
    \label{fig:ring_rate}
\end{figure}

We now perform a parameter sweep over increasing $L/m$.
The reference algorithms are gradient descent and triple momentum \citep{van2017fastest} with $\rho_{\text{gd}} = \frac{2}{m+L}$ and $\rho_{\text{TM}} = 1-\sqrt{\frac{m}{L}}$, are tuned to optimally converge  in the setting of no network dynamics. These algorithms can be described by
% \small
\begin{align}
    K_{\text{gd}}: \ & \mat{c|c}{
       I_c & -\frac{2}{m+L} I_c\hl
       I_c & 0 },  & 
    K_{\text{TM}}: \  &  
       \mat{cc|c}{
       I_c & (1+\frac{\rho_{\text{TM}} ^2}{2-\rho_{\text{TM}} })I_c & - \frac{1+\rho_{\text{TM}} }{L} I_c \\ 
       0 & \frac{\rho_{\text{TM}} ^2}{2-\rho_{\text{TM}} }I_c &  \frac{1+\rho_{\text{TM}} }{L}I_c \hl
       I_c & \frac{\rho_{\text{TM}} ^2}{(1+\rho_{\text{TM}})(2-\rho_{\text{TM}})}I_c & 0 }.
\end{align}
% \normalsize
% both of which are 
Figure \ref{fig:ring_rate} plots an upper-bound on the convergence rate \rw{$\rho$ as a function of $L/m$} developed by solving the analysis program in Theorem \ref{thm:alg_analysis} with filters \rw{ of length}  $\nu_{\max}=3$.
The compared curves in Figure \ref{fig:ring_rate} are  gradient descent, triple momentum, \rw{our proposed s}ynthesis  with identity filters in which $M_r$ may be different between modes $r \in \{1, \ldots, N_s\}$, and synthesis under an \rw{identity filter} with the additional conservatism-introducing constraint that the \rw{matrix $M_r$ in the LMI \eqref{eq:antipassivity_switched} is the same among all modes $r \in \{1, \ldots, N_s\}$}. Setting $M_r$ to the same value among all modes $r$ ignores the ring structure of the switching graph switching. All comparisons are performed with respect to $m=1$, in which $L$ is swept in the range $(1, 5]$. The \rw{synthesized} path-dependent controller \rw{has}  $\rho < 1$ $L < 3.360$, and the common storage controller is stable \rw{only} when $L < 1.660$. 
Gradient descent and triple momentum are both empirically unstable even when $L = 1.01, L' = 1.005.$ %\urg{confirm}.

\subsection{Time-Varying Delay: Increasing Delay}

Our next example involves synthesis of algorithms with increasing delay before the oracle $\nabla f$. \rw{We consider two network structures: Unrestricted Delay and  Packet Drops. Both networks involve a maximal delay $H_{\max}$ ranging from $H_{\max} = 0$ to $H_{\max} = 5$. The Unrestricted Delay setting involves a delay of $h_k \in \{0, \ldots, H_{\max}\}$, in which $h_k$ and $h_{k+1}$ are independent. 
The Packet Drop setting involves the delay increasing by 1 until the maximum delay of $H_{\max}$, or jumps to delay 0. This logic is visualized in Figure \ref{fig:logic:b} for $H_{\max}=3$.}

\begin{figure}[!h]
    \centering
    \includegraphics[width=0.9\linewidth]{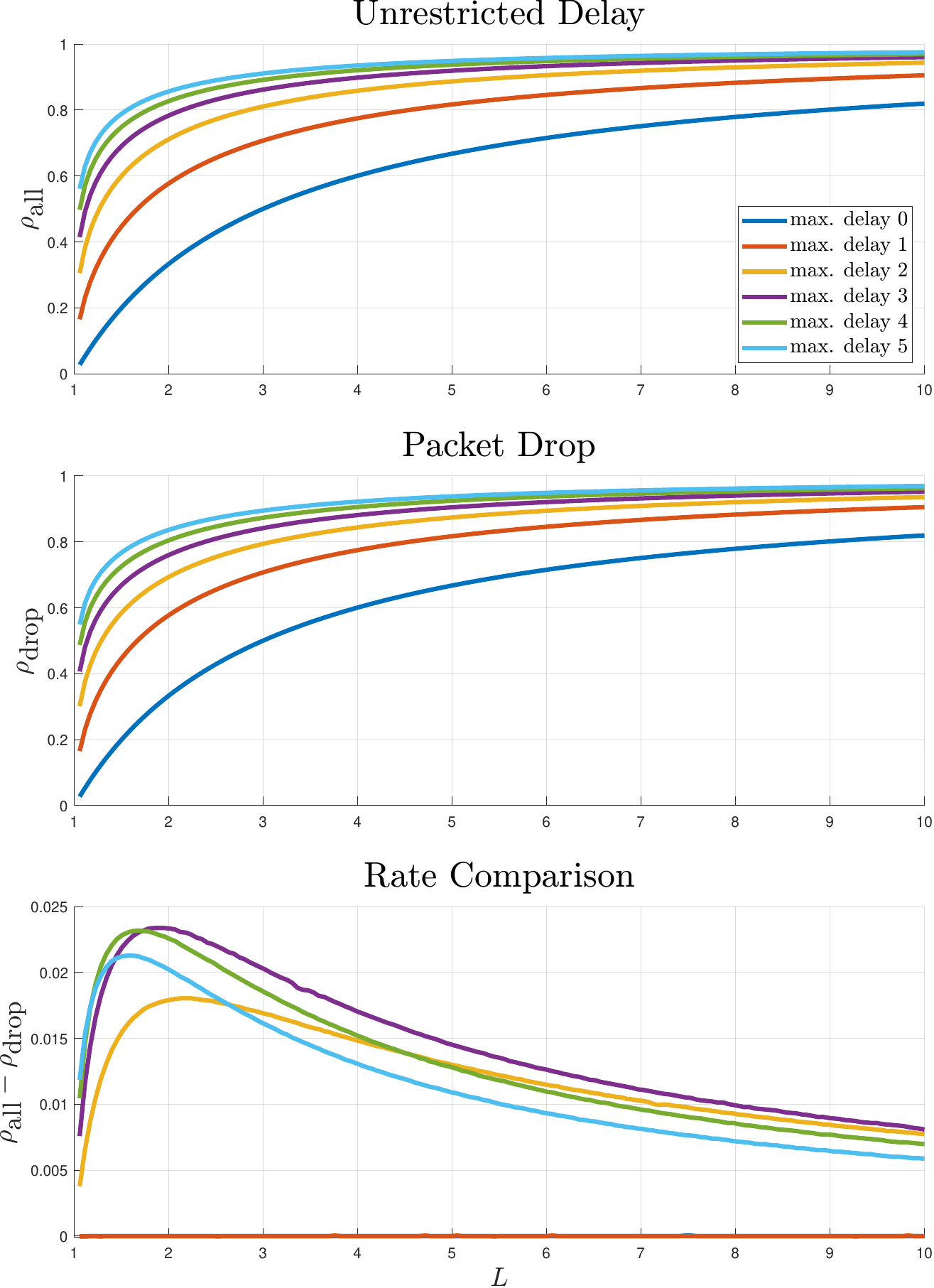}
    \caption{Exponential rates under increasing delay and  parameter $L$ }
    \label{fig:sweep_snap}
\end{figure}

Figure \ref{fig:sweep_snap} plots the computed convergence rate $\rho$ by synthesis under \rw{identity} multipliers $\lambda$ for $m = 1$ as $L$ is swept in the range $(1, 10]$. In each subplot, the color indicates the maximal delay $H_{\max}$ between $0$ and $5$. The top subplot shows the computed bounds on $\rho$ under unrestricted delay $(\rho_{\text{all}})$. The middle subplot depicts bounds on $\rho$ under packet drop logic $(\rho_{\text{drop}})$. The rates $\rho$ rise as the maximal delay increases. The bottom subplot draws the difference between the $(\rho_{\text{all}}) - (\rho_{\text{drop}})$. The computed rates of $H_{\max} \in \{0, 1\}$ are identical between the unrestricted delay and packet drop cases,\rw{when $H_{\max} \geq 2$ the packet-drop-synthesized controller has a smaller $\rho$ than a controller developed for the unrestricted delay case}. 
% When $H_{\max} > 2$, controllers synthesized for the packet drop offer a smaller convergence rate as compared to controllers synthesized for the }

% The maximal delay before the oracle increases from 0 to 6, and the upper-bound on the  convergence rate $\rho$ monotonically increases both in $L/m$ and in the delay. 

% \urg{
% Example content goes here.

% What would be sufficient?

% How much space will we have?

% }
\section{Conclusion}

\label{sec:conclusion}

\vspace{-0.25cm}

% An extended Arxiv version of this paper is available at \urg{[Arxiv link goes here] what will the arxiv version have that this work does not?}.

% This paper presented a delay-scheduled scheme to solve unconstrained optimization algorithms in the setting of time-varying discrete delays. The delay-scheduled optimization algorithm is posed as a (switched) family of linear systems, 
% each of which are in Lur\'e interconnection with the same gradient oracle (in class $\s_{m, L}$). 
This paper presented a method to analyze and synthesize switched optimization algorithms. 
Convergence to the optimal solution may be certified if the per-mode systems \rw{solve a regulator equation and achive stability (Theorem \ref{thm:convergence})}. Analysis occurs  through convex searches over Zames-Falb filter coefficients $\lambda$. Synthesis of switched algorithms given a network $P$ and fixed filters $\Psi(\lambda)$ is accomplished by using  regulator theory and output feedback techniques.
% Internal model control is applied to synthesize output feedback 
% With a fixed filter $\Psi(\lambda)$, output feedback synthesis is used to generate subcontrollers 
% With fixed multipliers, delay-schedules optimization algorithms can be synthesized through output feedback techniques. 
The analysis and synthesis tasks are jointly convex, and an alternating approach can be used to develop certified switching algorithms with exponential convergence guarantees.  
Future work involves the establishing of necessary conditions for \rw{convergence of  switched optimization algorithms}, and finding  delay-scheduled controllers in the  sampled data and unknown delay cases. Other work involves the formulating convex joint searches over controllers $R_r$ and filter coefficients $\lambda$, \rw{and extending our method synthesis distributed optimization algorithms.}

% The main generalization of this switched-system-based research involves the creation of optimization algorithms in the linear parameter varying setting. 
% \begin{ack}
% The authors would like to thank Zhiyu He, Chris Verhoek, Prof. Juan Francisco Camino, Prof. Roland T\'oth, Jaap Eising, and Alessio Moreschini for discussions about switched systems and optimization algorithms. 
% \end{ack}

\section*{Acknowledgements}
J. Miller and C. Scherer are by supported the Deutsche Forschungsgemeinschaft (DFG, German Research Foundation) under Germany's Excellence Strategy - EXC 2075 – 390740016. We acknowledge the support by the Stuttgart Center for Simulation Science (SimTech). F. Jakob acknowledges the support of the International Max Planck Research School for Intelligent Systems (IMPRS-IS).

\bibliographystyle{siam}
\bibliography{references}

@book{ben2001lectures,
  title={Lectures on modern convex optimization: analysis, algorithms, and engineering applications},
  author={Ben-Tal, Aharon and Nemirovski, Arkadi},
  year={2001},
  publisher={SIAM}
}

@incollection{vandenberghe2005interior,
  title={Interior-point algorithms for semidefinite programming problems derived from the {KYP} lemma},
  author={Vandenberghe, Lieven and Balakrishnan, V Ragu and Wallin, Ragnar and Hansson, Anders and Roh, Tae},
  booktitle={Positive polynomials in control},
  pages={195--238},
  year={2005},
  publisher={Springer}
}

@article{scherer2025tutorial,
  title={{A Tutorial on Convex Design of Optimization Algorithms by Integral Quadratic Constraints}},
  author={Scherer, Carsten W and Ebenbauer, Christian},
  journal={Annual Review of Control, Robotics, and Autonomous Systems},
  volume={9},
  year={2025},
  publisher={Annual Reviews}
}

@article{rotaru2024exact,
  title={Exact worst-case convergence rates of gradient descent: a complete analysis for all constant stepsizes over nonconvex and convex functions},
  author={Rotaru, Teodor and Glineur, Fran{\c{c}}ois and Patrinos, Panagiotis},
  journal={arXiv preprint arXiv:2406.17506},
  year={2024}
}

@book{zhou1998essentials,
  title={Essentials of robust control},
  author={Zhou, Kemin and Doyle, John Comstock},
  volume={104},
  year={1998},
  publisher={Prentice hall Upper Saddle River, NJ}
}

@article{zattoni2013output,
  title={The output regulation problem with stability for linear switching systems: A geometric approach},
  author={Zattoni, Elena and Perdon, Anna Maria and Conte, Giuseppe},
  journal={Automatica},
  volume={49},
  number={10},
  pages={2953--2962},
  year={2013},
  publisher={Elsevier}
}

@article{lur1944theory,
  title={On the theory of stability of control systems},
  author={Lur’e, Anatoliy Isakovich and Postnikov, Vladimir N},
  journal={Applied mathematics and mechanics},
  volume={8},
  number={3},
  pages={246--248},
  year={1944}
}

@article{Ahmadi_2014,
   title={{Joint Spectral Radius and Path-Complete Graph Lyapunov Functions}},
   volume={52},
   ISSN={1095-7138},
   number={1},
   journal={SIAM J. Control Optim.},
   publisher={Society for Industrial & Applied Mathematics (SIAM)},
   author={Ahmadi, Amir Ali and Jungers, Rapha\"el M. and Parrilo, Pablo A. and Roozbehani, Mardavij},
   year={2014},
   month=jan, pages={687–717} }

@article{apkarian1995self,
  title={Self-scheduled {H}-infinity control of linear parameter-varying systems: a design example},
  author={Apkarian, Pierre and Gahinet, Pascal and Becker, Greg},
  journal={Automatica},
  volume={31},
  number={9},
  pages={1251--1261},
  year={1995},
  publisher={Elsevier}
}

@article{SCHERER2023robustexp,
title = {{Robust Exponential Stability and Invariance Guarantees with General Dynamic O'Shea-Zames-Falb Multipliers}},
journal = {IFAC-PapersOnLine},
volume = {56},
number = {2},
pages = {5799-5804},
year = {2023},
note = {22nd IFAC World Congress},
author = {Carsten W. Scherer},
}

@inproceedings{conte2020modeling,
  title={Modeling discrete time systems with variable delays as switching systems without delays},
  author={Conte, Giuseppe and Perdon, Anna Maria and Zattoni, Elena},
  booktitle={2020 European Control Conference (ECC)},
  pages={1591--1594},
  year={2020},
  organization={IEEE}
}

@article{de2012gain,
  title={Gain-scheduled dynamic output feedback control for discrete-time {LPV} systems},
  author={De Caigny, Jan and Camino, Juan F and Oliveira, Ricardo CLF and Peres, Pedro LD and Swevers, Jan},
  journal={Int. J. Robust Nonlinear Control.},
  volume={22},
  number={5},
  pages={535--558},
  year={2012},
  publisher={Wiley Online Library}
}

@INPROCEEDINGS{gahinet1993lmilab,
  author={Gahinet, P. and Nemirovskii, A.},
  booktitle={	Proc. IEEE Conf. Decis. Control}, 
  title={{General-Purpose {LMI} solvers with benchmarks}}, 
  year={1993},
  volume={},
  number={},
  pages={3162-3165 vol.4}}

@article{ramaswamy2021optimization,
  title={Optimization over time-varying networks and unbounded information delays},
  author={Ramaswamy, Arunselvan and Redder, Adrian and Quevedo, Daniel E},
  journal={IEEE Transactions on Automatic Control},
  volume={67},
  number={8},
  pages={4131--4137},
  year={2021},
  publisher={IEEE}
}

@article{doostmohammadian2021consensus,
  title={Consensus-based distributed estimation in the presence of heterogeneous, time-invariant delays},
  author={Doostmohammadian, Mohammadreza and Khan, Usman A and Pirani, Mohammad and Charalambous, Themistoklis},
  journal={IEEE Control Syst. Lett.},
  volume={6},
  pages={1598--1603},
  year={2021},
  publisher={IEEE}
}

@article{stoorvogel2000performance,
  title={Performance with regulation constraints},
  author={Stoorvogel, Anton A and Saberi, Ali and Sannuti, Peddapullaiah},
  journal={Automatica},
  volume={36},
  number={10},
  pages={1443--1456},
  year={2000},
  publisher={Elsevier}
}

@article{briat2009delay,
  title={Delay-scheduled state-feedback design for time-delay systems with time-varying delays—a {LPV} approach},
  author={Briat, Corentin and Sename, Olivier and Lafay, Jean-Fran{\c{c}}ois},
  journal={Systems \& Control Letters},
  volume={58},
  number={9},
  pages={664--671},
  year={2009},
  publisher={Elsevier}
}

@book{desoer1975feedback,
  title={Feedback Systems: Input-Output Approach},
  author={Desoer, Charles A and Vidyasagar, Mathukumalli},
  year={1975},
  publisher={Academic Press, London}
}

@article{francis1976internal,
  title={The internal model principle of control theory},
  author={Francis, Bruce A and Wonham, Walter Murray},
  journal={Automatica},
  volume={12},
  number={5},
  pages={457--465},
  year={1976},
  publisher={Elsevier}
}

@article{conte2021disturbance,
  title={Disturbance decoupling and model matching problems for discrete-time systems with time-varying delays},
  author={Conte, G and Perdon, AM and Zattoni, E and Animobono, D},
  journal={Nonlinear Analysis: Hybrid Systems},
  volume={41},
  pages={101043},
  year={2021},
  publisher={Elsevier}
}

@article{korouglu2008lpv,
  title={{LPV control for robust attenuation of non-stationary sinusoidal disturbances with measurable frequencies}},
  author={K{\"o}ro{\u{g}}lu, Hakan and Scherer, Carsten W},
  journal={IFAC Proceedings Volumes},
  volume={41},
  number={2},
  pages={4928--4933},
  year={2008},
  publisher={Elsevier}
}

@article{lessard2016analysis,
  title={Analysis and design of optimization algorithms via integral quadratic constraints},
  author={Lessard, Laurent and Recht, Benjamin and Packard, Andrew},
  journal={SIAM J. Optim.},
  volume={26},
  number={1},
  pages={57--95},
  year={2016},
  publisher={SIAM}
}

@inproceedings{scherer2023optimization,
  title={Optimization algorithm synthesis based on integral quadratic constraints: A tutorial},
  author={Scherer, Carsten W. and Ebenbauer, Christian and Holicki, Tobias},
  booktitle={	Proc. IEEE Conf. Decis. Control},
  pages={2995--3002},
  year={2023},
  organization={IEEE}
}

@article{van2017fastest,
  title={The fastest known globally convergent first-order method for minimizing strongly convex functions},
  author={Van Scoy, Bryan and Freeman, Randy A and Lynch, Kevin M},
  journal={IEEE Control Systems Letters},
  volume={2},
  number={1},
  pages={49--54},
  year={2017},
  publisher={IEEE}
}

@article{jakob2025linear,
  title={{A Linear Parameter-Varying Framework for the Analysis of Time-Varying Optimization Algorithms}},
  author={Jakob, Fabian and Iannelli, Andrea},
  journal={arXiv:2501.07461},
  year={2025}
}

@article{scherer2021convex,
  title={Convex synthesis of accelerated gradient algorithms},
  author={Scherer, Carsten and Ebenbauer, Christian},
  journal={SIAM J. Control Optim.},
  volume={59},
  number={6},
  pages={4615--4645},
  year={2021},
  publisher={SIAM}
}

@inproceedings{holicki2021algorithm,
  title={Algorithm design and extremum control: Convex synthesis due to plant multiplier commutation},
  author={Holicki, Tobias and Scherer, Carsten W},
  booktitle={	Proc. IEEE Conf. Decis. Control},
  pages={3249--3256},
  year={2021},
  organization={IEEE}
}

@ARTICLE{jovanovic2021robustness,
  author={Mohammadi, Hesameddin and Razaviyayn, Meisam and Jovanović, Mihailo R.},
  journal={IEEE Transactions on Automatic Control}, 
  title={{Robustness of Accelerated First-Order Algorithms for Strongly Convex Optimization Problems}}, 
  year={2021}
}

@article{michalowski2021,
author = {Simon Michalowsky and Carsten Scherer and Christian Ebenbauer},
title = {Robust and structure exploiting optimisation algorithms: an integral quadratic constraint approach},
journal = {International Journal of Control},
volume = {94},
number = {11},
pages = {2956--2979},
year = {2021},
publisher = {Taylor \& Francis}}

@INPROCEEDINGS{wen2008switched,
  author={Zhang Wen'an and Yu Li and Song Hongbo},
  booktitle={2008 27th Chinese Control Conference}, 
  title={A switched system approach to networked control systems with time-varying delays}, 
  year={2008}}

@article{Megretski1997,
    author={A. Megretski and A. Rantzer},
    title={{System analysis via Integral Quadratic Constraints}},
    journal={IEEE Transactions on Automatic Control},
    year={1997},
}

@inproceedings{BoczarLessardRecht,
    author={R. Boczar and L. Lessard and B. Recht},
    title={{Exponential Convergence Bounds Using Integral Quadratic Constraints}},
    booktitle={American Control Conference},
    year={2018},
}

@Article{Daafouz2002,
  author    = {Daafouz, J. and Riedinger, P. and Iung, C.},
  journal   = {IEEE Transactions on Automatic Control},
  title     = {Stability analysis and control synthesis for switched systems: a switched Lyapunov function approach},
  year      = {2002},
  issn      = {0018-9286},
  month     = nov,
  number    = {11},
  pages     = {1883--1887},
  volume    = {47},  
  publisher = {Institute of Electrical and Electronics Engineers (IEEE)},
}

@Article{Lee2020,
  author    = {Lee, Donghwan and Dullerud, Geir E. and Hu, Jianghai},
  journal   = {Automatica},
  title     = {{Graph Lyapunov function for switching stabilization and distributed computation}},
  year      = {2020},
  issn      = {0005-1098},
  month     = jun,
  pages     = {108923},
  volume    = {116},
  publisher = {Elsevier BV},
}

@Article{Zhao2020,
  author    = {Zhao, Ying and Liu, Yu and Ma, Dan},
  journal   = {IEEE Transactions on Circuits and Systems I: Regular Papers},
  title     = {Output Regulation for Switched Systems With Multiple Disturbances},
  year      = {2020},
  issn      = {1558-0806},
  month     = dec,
  number    = {12},
  pages     = {5326--5335},
  volume    = {67},
  
  publisher = {Institute of Electrical and Electronics Engineers (IEEE)},
}

@Article{Lee2006,
  author    = {Lee, Ji-Woong and Dullerud, Geir E.},
  journal   = {Automatica},
  title     = {Uniform stabilization of discrete-time switched and Markovian jump linear systems},
  year      = {2006},
  issn      = {0005-1098},
  month     = feb,
  number    = {2},
  pages     = {205--218},
  volume    = {42},

  publisher = {Elsevier BV},
}

\end{document}